\documentclass[10pt]{article}
\usepackage{amssymb,amsfonts,amsmath, amsthm, amsbsy}
\usepackage{color,graphicx, subcaption, placeins, array, mathtools, url, color, rotating} 
\usepackage[text={6.75in,9.5in},centering,letterpaper]{geometry}
\usepackage[hypertexnames=false,colorlinks=true,urlcolor=blue,linkcolor=blue,citecolor=blue]{hyperref}
\usepackage{tikz}
\usepackage[shortlabels]{enumitem}
\usepackage[margin = 10pt, font=small, labelfont=bf, labelsep = endash]{caption}
\usepackage{comment}

\everymath={\displaystyle}
 \numberwithin{equation}{section}

\newtheorem{Theorem}{Theorem}[section]

\newtheorem{Remark}[Theorem]{Remark}

\usepackage{todonotes}
\newcommand{\PvHx}[1]{\todo[color=green]{PvH: #1}} 
\newcommand{\PvH}[1]{\todo[inline,color=green]{PvH: #1}}

\newcommand{\eps}{\varepsilon}

\title{Deformations of acid-mediated invasive  tumors in a model with Allee effect}
\author{Paul Carter\thanks{Department of Mathematics, University of California, Irvine, USA}\and Arjen Doelman\thanks{Mathematical Institute, Leiden University, Leiden, Netherlands}\and Peter van Heijster\thanks{Mathematical and Statistical Methods—Biometris, Wageningen University \& Research, Wageningen, Netherlands}\and Daniel Levy\thanks{Program in Applied and Computational Mathematics, Princeton University, Princeton, USA }\and Philip Maini\thanks{Mathematical Institute, University of Oxford, Oxford, UK}\and Erin Okey\thanks{Division of Applied Mathematics, Brown University, Providence, RI}\and  Paige Yeung\thanks{Massachusetts Institute of Technology, Cambridge, USA. }}

\begin{document}

\maketitle

\begin{abstract}
We consider a Gatenby--Gawlinski-type model of invasive tumors in the presence of an Allee effect. We describe the construction of bistable one-dimensional traveling fronts using singular perturbation techniques in different parameter regimes corresponding to tumor interfaces with, or without, acellular gap. By extending the front as a planar interface, we perform a stability analysis to long wavelength perturbations transverse to the direction of front propagation and derive a simple stability criterion for the front in two spatial dimensions. In particular we find that in general the presence of the acellular gap indicates transversal instability of the associated planar front, which can lead to complex interfacial dynamics such as the development of finger-like protrusions and/or different invasion speeds.
\end{abstract}








\section{Introduction}
Evolving cancerous tumors go through several phases. Tumors initially appear as relatively small congregations of cells that have undergone genetic mutations that typically causes the tumor to grow at an abnormally high rate~\cite{wodarz2014dynamics}. At onset, this growth is regular and sustained by nutrients, etc., that are transported to the tumor by the governing local diffusion processes. However, at a certain stage, this regular mechanism is no longer able to drive the tumor growth process and the tumor enters the next phase in which it invades the surrounding, healthy, tissue. This stage is characterized by a deformation of the surface of the original clump of (tumorous) cells. The nature of resulting morphology can determine the invasive capacity, and therefore the overall severity, of the tumor. To invade tissue, tumor cells need to overcome the physical barriers of normal cells and extracellular matrix (ECM), the matrix that provides a physical scaffolding for cells. Understanding how tumor cells achieve this is of obvious importance, especially since at the next stage of the process -- known as metastasis -- the invading surface of the preceding stage may break away from the primary tumor and invade other parts of body to form, often fatal, secondary tumors. 

There are many ways to model the multi-scale process of tumor growth and invasion mathematically, ranging from fully discrete to  continuum models, and various hybrid models in between -- see for instance \cite{anderson2008integrative,lowengrub2010nonlinear,turner2002intercellular}. Here, we take the continuum point of view and focus on (generalized, nonlinear) reaction-diffusion type models. Several of these models address the role of matrix metalloproteinases (MMPs), which are enzymes secreted by tumor cells that degrade ECM. For example, Perumpanani {\em et al.} (1996) \cite{Perumpanani96} proposed a model consisting of six coupled partial differential equations (PDEs) for cell, ECM and MMP densities, where there are cells of different types (phenotypes) and movement is via the process of haptotaxis (movement up adhesive gradients) as well as diffusion. They analysed, in one spatial dimension, the traveling wave behavior of this system, exploring how the wave speed depends on various cell properties. More recently, Katsaounis {\em et al.} (2024) \cite{Katsaounis24} developed a hybrid multiscale 3D model employing PDEs and stochastic differential equations in which cell diffusion is modelled nonlinearly. They carry out a numerical study of the system in the context of a number of biologically motivated case studies. In 1996, Gatenby and Gawlinski \cite{GG} published their seminal paper in which they stated their acid-mediated invasion hypothesis. Namely, they proposed that tumor cells, undergoing glycolysis, which produces lactic acid, could invade normal cells due to the acid being more toxic to normal cells than to tumor cells. Their PDE model was composed of three coupled equations, with a cross-dependent degenerate diffusion term for tumor cell movement. This led to the formation of sharp interfaces representing the surface of the growing tumor and the prediction that, in certain cases, an acellular gap would appear between the tumor and normal cell populations. This prediction was validated in the context of head and neck cancer. This model has been extended in a number of studies to account also for ECM degradation (see, for example, \cite{Martin10,Strobl20}).

\begin{figure}
\centering
\begin{minipage}{0.3\textwidth}
\includegraphics[width=1\linewidth]{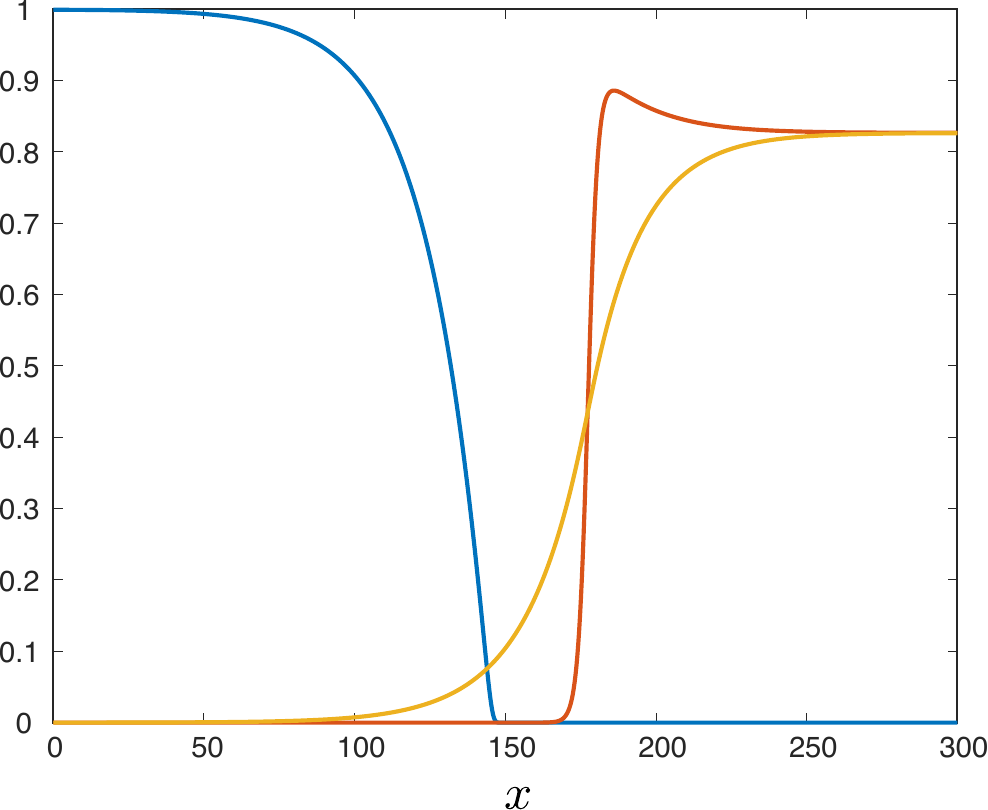}
\end{minipage}
\hspace{0.05\textwidth}
\begin{minipage}{0.6\textwidth}
\includegraphics[width=0.325\linewidth]{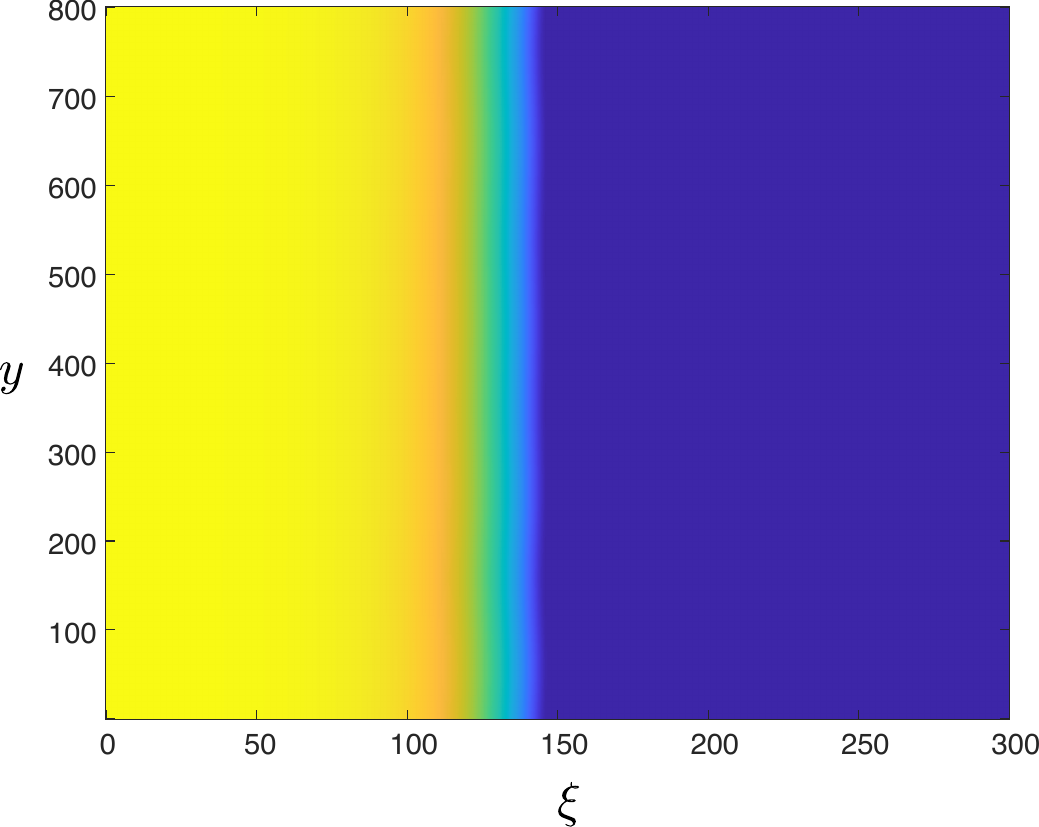}
\includegraphics[width=0.325\linewidth]{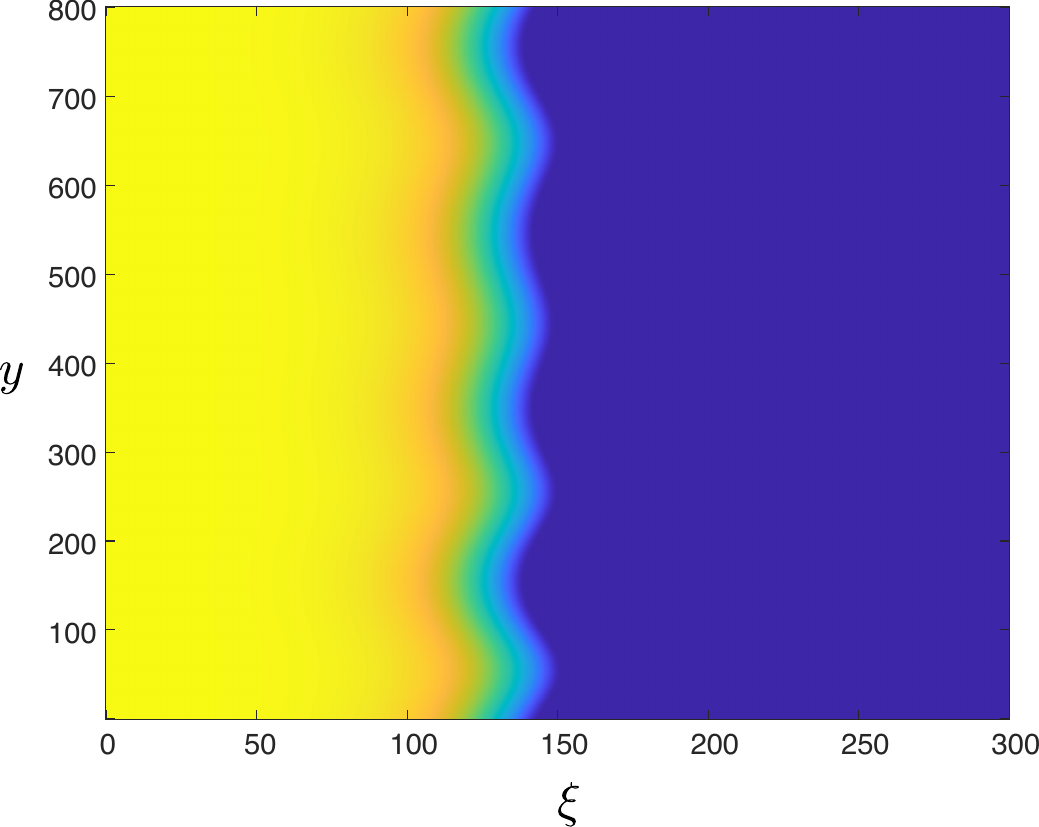}
\includegraphics[width=0.325\linewidth]{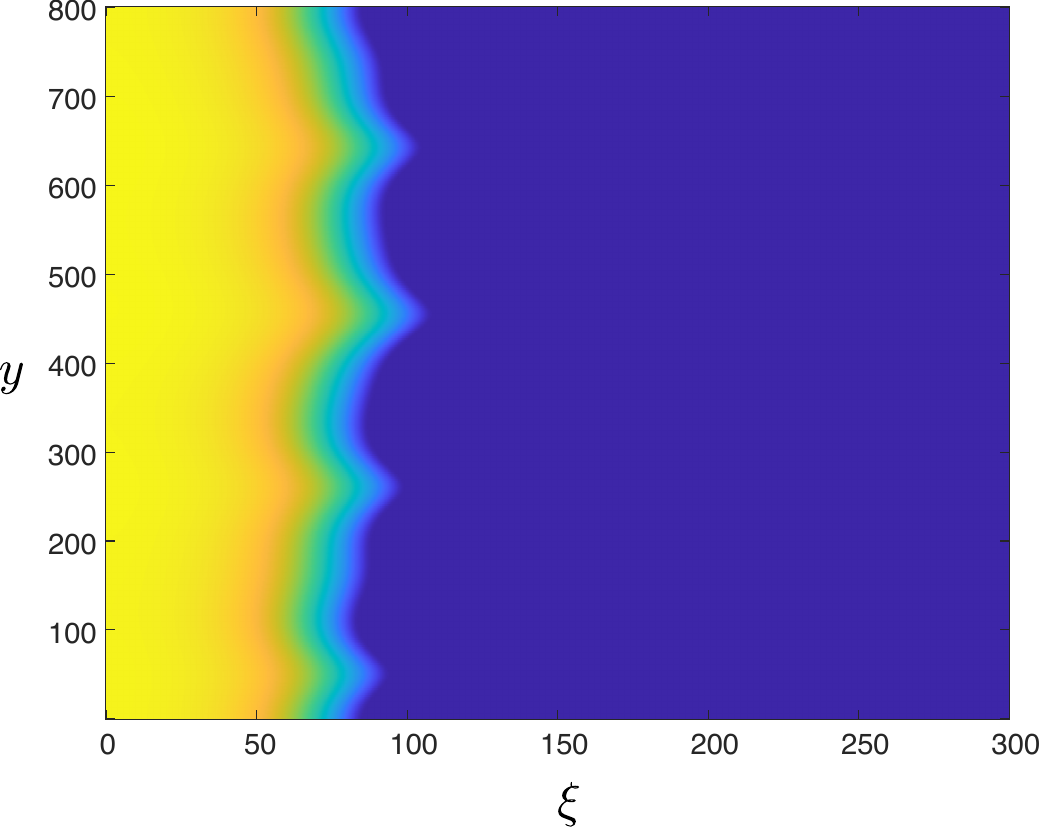}\\
\includegraphics[width=0.325\linewidth]{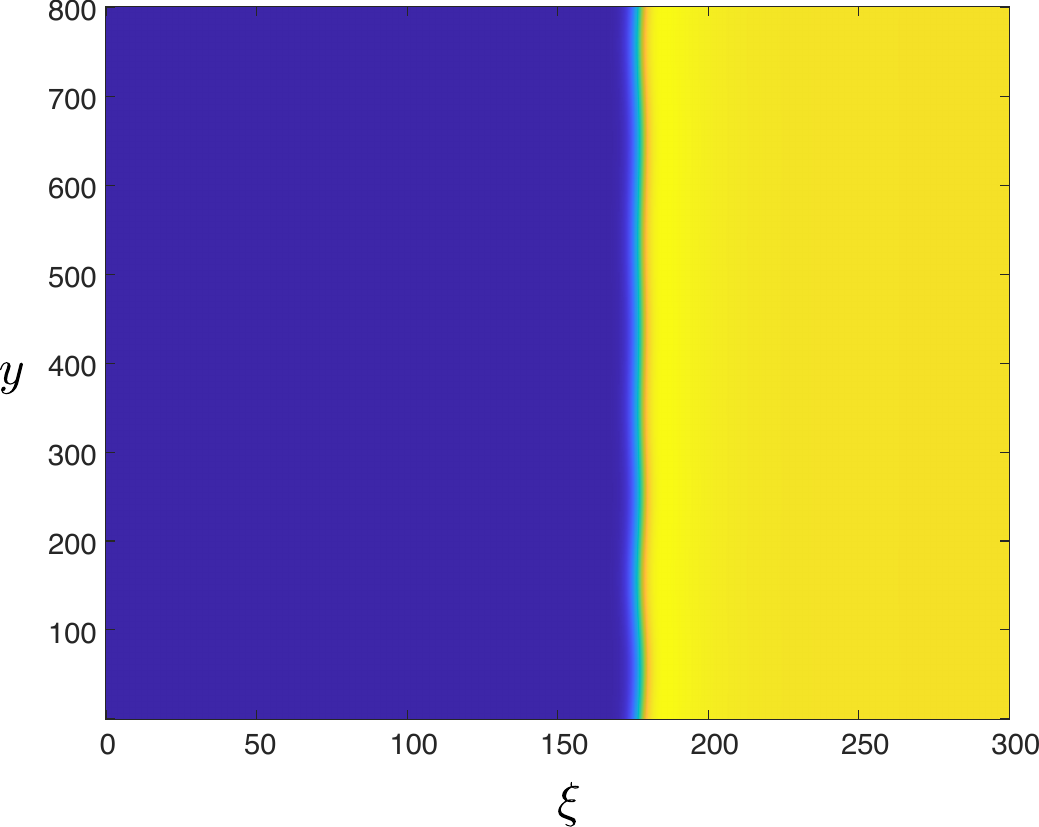}
\includegraphics[width=0.325\linewidth]{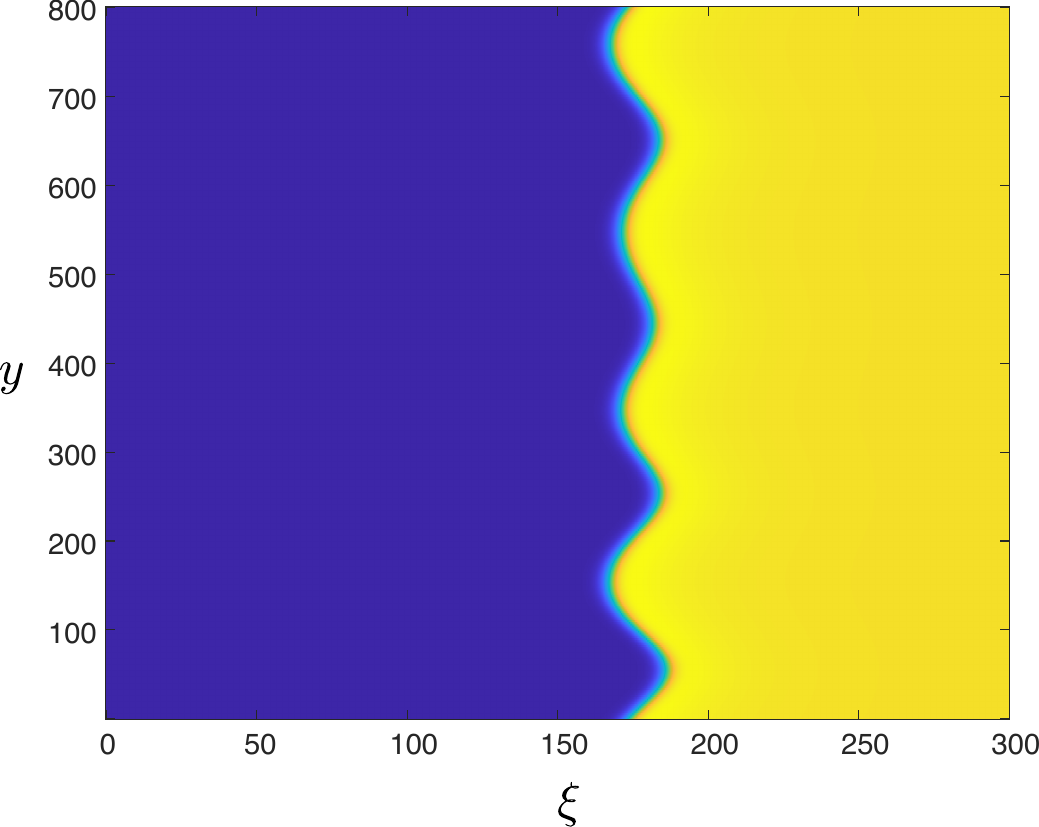}
\includegraphics[width=0.325\linewidth]{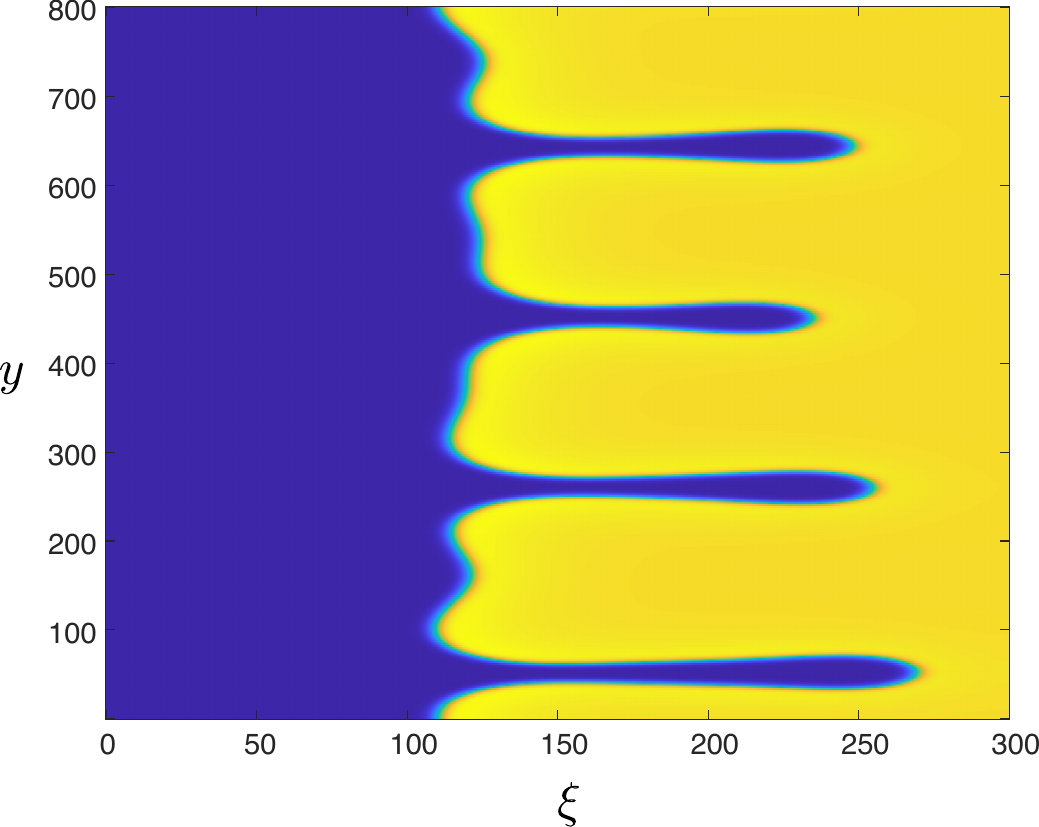}
\end{minipage}
\caption{Results of a direct numerical simulation of~\eqref{eq:GG} exhibiting the transversal long wavelength instability for parameter values taken from~\cite{GG}  $(a,\kappa, \delta_1, \delta_2, \delta_3, \rho, \eps) = (0.35,0.1,12.5,0.1,70.0,1.0, 0.0063)$. In the left panel, the underlying (longitudinally stable) one-dimensional traveling wave is shown (in which the normal cell density $U$ is plotted in blue, the tumor cell density $V$ in red, and the acid concentration $W$ in yellow. The six right panels depict the results of two-dimensional simulations. The initial conditions of the two-dimensional run were constructed by trivially extending this $1$D-stable profile in the $y$-direction and adding a small amount of noise. The simulations were performed in a co-moving frame $\xi=x+ct$; in the laboratory frame the tumor would travel with speed $c=0.0401$ to the left, corresponding to the wave speed of the initial front profile (so that in the absence of any instability, the fronts would appear stationary). The corresponding tumor $U$-profiles (top three panels) and $V$-profiles (bottom three panels) are plotted (for $t=20000,30000,40000$ from left to right); for all profiles yellow indicates high density of cells, and blue depicts low density of cells. }
\label{f:Intro-interfaces}
\end{figure}

The mathematical challenges raised by the nonlinear cross-dependent diffusion  term in \cite{GG} have inspired many mathematical studies of the traveling wave behavior of this system in one spatial dimension (see \cite{davis2022traveling, Fasano09, Gallay22} and the references therein). However, these in essence one-dimensional interfaces can only be expected to be stable -- and thus observable -- in the initial `regular' stage of tumor growth. In the follow-up phase in which the tumor starts to invade the surrounding tissue, the interface most likely will evolve into a structure with a two- or three-dimensional nature. Viewing such shapes as possibly arising from a spatial patterning instability, but departing from the idea of a (growing) tumor surface governed by traveling waves, Chaplain {\em et al.} (2001) \cite{Chaplain01} proposed a Turing-type model composed of growth activating and growth inhibiting chemicals and showed how these gave rise to spatial patterns on a spherical surface which, they proposed, would induce `columnar outgrowths of invading cancer cells'. In this paper, we propose a different mechanism for the formation of outgrowths of invading tumor cells that is based on the tumor surface as the interface between tumorous and healthy tissue governed by bistable traveling fronts (see Fig. \ref{f:Intro-interfaces}). Namely, we analyse front evolution -- or interface dynamics -- in a slightly modified form of the Gatenby-Gawlinski \cite{GG} model to explore the potential appearance of transversal instabilities on interfaces that are longitudinally stable. In other words, we consider two-dimensional fronts for which the underlying one-dimensional traveling waves are stable in the direction of propagation, but for which instabilities may form in directions transversal to this direction. In Fig. \ref{f:Intro-interfaces} we show the outcome of a simulation of the (slightly modified) Gatenby-Gawlinski model in two space dimensions that indeed exhibits such a transversal instability: the interface develops `fingers' similar to the phenomenon of viscous fingering in fluid dynamics. We emphasize that the instability develops purely due to extending the (longitudinally stable) front as an interface in two spatial dimensions, in the absence of any deterministic or stochastic forcing or inhomogeneity.

In our modified version of the Gatenby-Gawlinski model \cite{GG}, we assume that tumor cells are not fully resistant to lactic acid by including a death term in the tumor cell equation. We also regularise the nonlinear diffusion term by assuming that normal cells do not form an impenetrable barrier to tumor cell invasion. While we do this primarily for mathematical convenience, it should be noted that there is a biological justification for this as some tumor cells upregulate the expression of aquaporins (see the review \cite{Verkman08}) which, it has been hypothesised, may allow cells to `bulldoze' their way through surrounding cells and tissues \cite{McLennan20}. Lastly, we replace the logistic growth term for cancer cells by a term that allows for the Allee effect, which has been proposed to play a role in modeling cancer tumors similar to modeling ecosystems -- see \cite{Gerlee22, Johnson19, korolev2014turning, Sewalt16} and the references therein. Thus, we consider the following reaction-diffusion model for tumor invasion,
\begin{align}
\begin{split}
\label{eq:GG}
U_\tau &=U(1-U)-\delta_1 UW\,,\\
V_\tau  &= \rho V(1-V)(V-a)-\delta_2 VW+\nabla \cdot ((1+\kappa -U) \nabla V)\,,\\
W_\tau &= \delta_3(V-W)+\frac{1}{\eps^2}\Delta W\,,
\end{split}
\end{align}
where $(U,V,W)(x,y,\tau)$ represent normal cell density, tumor cell density, and acid (concentration) produced by the tumor cells, respectively, at spatial position $(x,y)\in\mathbb{R}^2$ and time $\tau\in \mathbb{R}^+$. The parameter $\delta_1$ measures the destructive influence of acid on healthy tissue and can thus be seen as an indicator of tumor aggressivity, the new parameter $\delta_2$ represents the impact of acid on the tumor itself and we assume that $0\leq \delta_2 <\delta_1$, i.e. the effect of $W$ on $V$ is less than that of $W$ on $U$ (but not necessarily zero). The parameter $\rho > 0$ measures the relative production rate of tumor cells compared to healthy cells and $0<a<1$ the relative strength of the Allee effect. Since the diffusive spreading speed of acid is much faster than that of cells \cite{GG} it follows that $0< \eps \ll 1$ so that reaction-diffusion system (\ref{eq:GG}) is singularly perturbed. Finally, the parameter $\kappa$ measures the regularized nonlinear diffusion effect and is assumed to satisfy $0<\eps\ll \kappa\ll1$.

In the spirit of \cite{korolev2014turning}, the present paper was motivated by recent progress on the formation of invasive, `fingering', interface patterns as in Fig. \ref{f:Intro-interfaces} in the context of coexistence patterns (between grasslands and bare soil) in dryland ecosystems \cite{carter24, CDLOR, fernandez2019front}. In fact, in \cite{carter24,CDLOR} criteria have been developed by which the transversal (in)stability of planar interfaces in a general class of singularly perturbed two-component reaction–diffusion equations -- that includes typical dryland ecosystem models -- can be determined. The present model does not directly fall into this class of systems: \eqref{eq:GG} contains three components and has a nonlinear diffusion term -- unlike the systems considered in \cite{CDLOR}. Nevertheless, we show in this paper that the methods by which the criteria in \cite{CDLOR} are deduced can also be applied to the present modified Gatenby-Gawlinski model.    

\begin{figure}[t]
\centering
\includegraphics[width=0.4\linewidth]{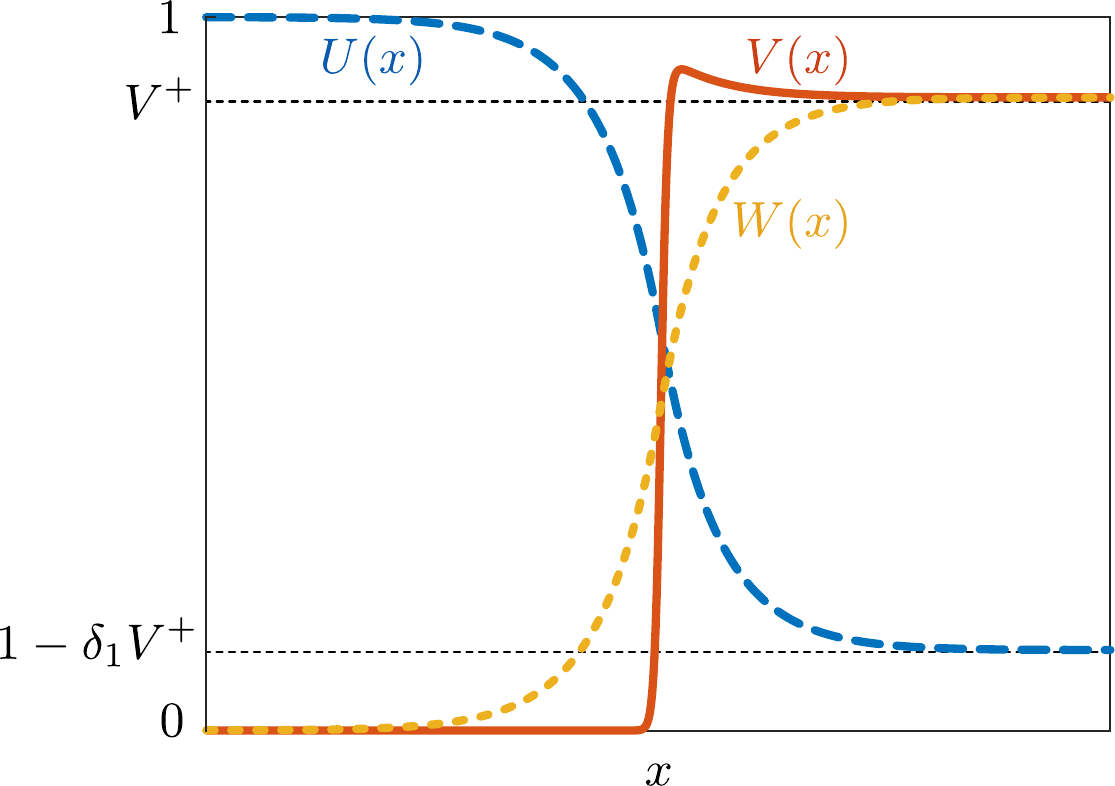}\hspace{0.1\linewidth}
\includegraphics[width=0.37\linewidth]{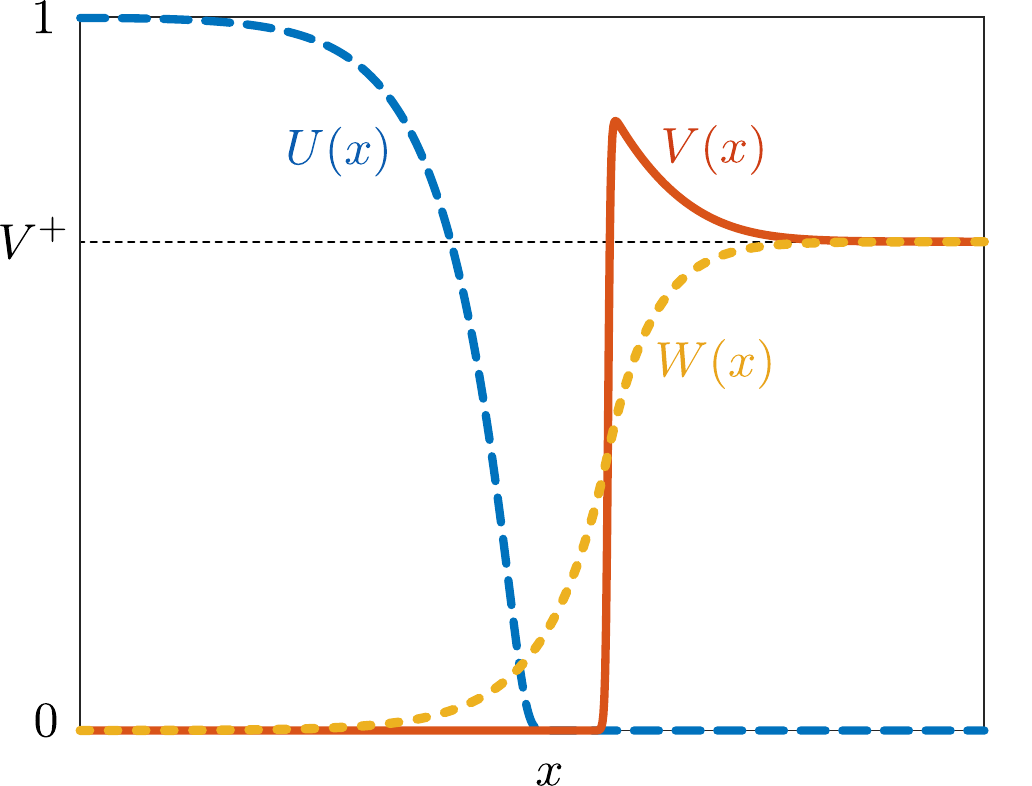}
\caption{Bistable traveling front solutions in the benign (left) and malignant (right) cases. Note that in the malignant case, an acellular gap may appear in between the interfaces formed by the $U$ and $V$ profiles and that this case corresponds to the situation of Fig. \ref{f:Intro-interfaces}. 
}
\label{f:front}
\end{figure}

As in \cite{CDLOR}, the transversal (in)stability results are derived from a detailed investigation into the structure of the (one-dimensional) traveling front on which the evolving two-dimensional interface is based (see Fig. \ref{f:Intro-interfaces}). These fronts correspond to heteroclinic connections between  stable background states -- or homogeneous equilibria -- in \eqref{eq:GG}. The background states of \eqref{eq:GG} are given by
\begin{align}\label{eq:steadystates}
P_1 = (0,0,0), \quad P_2 = (1,0,0), \quad P_3^\pm = (1-\delta_1 V^\pm, V^\pm, V^\pm), \quad P_4^\pm=(0,V^\pm, V^\pm),
\end{align}
with
\begin{align} \label{Vplusmin}
V^\pm = \frac{\rho(1+a)-\delta_2 \pm\sqrt{\left(\rho(1+a)-\delta_2\right)^2-4\rho^2a}}{2\rho},
\end{align}
where we note that the states $P_3^\pm$ only are biologically relevant for parameter combinations such that $V^\pm$ are real and positive and $\delta_1V^+<1$, and $P_4^\pm$ for parameters such that $V^\pm > 0$. It follows from linear stability analysis (see Appendix~\ref{app:steadystates}) that $P_2$ is stable, while $P_1, P_3^-, P_4^-$ are unstable. The background state $P_3^+$ is stable for $\delta_1V^+<1$, while $P_4^+$ is stable if $\delta_1V^+>1$. In this paper, we refer to the former as the \emph{benign} case, and the latter as the \emph{malignant} case, based on whether normal cells $U$ can coexist with tumor cells $V$ in the nontrivial stable steady state. In the benign case $\delta_1V^+<1$ we will construct (and further study) a bistable traveling front between homogeneous background states $P_2$ and $P_3^+$, while in the malignant case $\delta_1V^+>1$, we will consider a bistable front between $P_2$ and $P_4^+$ -- see Fig.~\ref{f:front}.

We construct three types of bistable traveling fronts/tumor interfaces using geometric singular perturbation theory: a benign front, a malignant no-gap front and a malignant gap front (see upcoming Fig. \ref{f:3D_singular} for the distinct geometries of these cases). The overall geometry of the construction of fronts is similar to that of~\cite{davis2022traveling}, in which invasion fronts were constructed in the model~\eqref{eq:GG} in the absence of the Allee effect. Some differences arise here due to the presence of the Allee effect, which results in bistable fronts existing at a unique wave speed, as opposed to a range of speeds. These bistable fronts are more amenable to a two-dimensional stability analysis, as the critical spectrum associated with such fronts in one spatial dimension takes the form of a single eigenvalue at $\lambda=0$ due to translation invariance. As in \cite{CDLOR}, we do not explicitly analyze the longitudinal, i.e. one-dimensional (in the direction of propagation), stability of these traveling fronts: we assume that they are longitudinally stable and thus that the translational eigenvalue at $\lambda = 0$ is the most critical eigenvalue. Note that the assumptions are based on numerical observations and especially on numerical evaluations of the spectrum associated with the linearized stability problem (cf. \S\ref{S:NUM}). The stability of the associated planar interface spanned by the one-dimensional fronts can be determined by an additional transversal Fourier expansion -- parameterized by wavenumber $\ell \in \mathbb{R}$. As a consequence, the translational eigenvalue re-appears as a local extremum at $\ell = 0$ of a critical curve $\lambda_\mathrm{c}(\ell)$, 
i.e. $\lambda_\mathrm{c}(0) = 0$. Observe that this curve is symmetric in $\ell$. Extending the methods developed in \cite{CDLOR}, we derive an explicit approximation of $\lambda_\mathrm{c}''(0)$, where prime denotes differentiation with respect to $\ell$. Clearly, the interface is unstable with respect to transversal long wavelength perturbations if $\lambda_\mathrm{c}''(0) > 0$ (so that $\lambda_\mathrm{c}(\ell) > 0$ for $\ell$ sufficiently small). In that case, the originally planar interface will be unstable and growing long wavelength spatial structures will emerge from the interface. Numerical simulations show that the interface may develop protruding fingers (Fig. \ref{f:Intro-interfaces}) or cusps (see~\S\ref{S:NUM}), which is similar to the observations of evolving interfaces (between grasslands and bare soil) in dryland ecosystems under the same circumstances (i.e. also with $\lambda_\mathrm{c}''(0) > 0$) \cite{carter24, CDLOR,fernandez2019front}.   

The explicit character of our analysis establishes a direct correspondence between the nature -- benign, malignant no-gap/gap -- of the underlying traveling front and the initiation of transversal instabilities of the evolving interface between tumorous and regular cells. We briefly summarize some of our main findings (and refer to~\S\ref{S:NUM} for the more extensive discussion): Firstly, we deduce directly from our derivation of $\lambda_\mathrm{c}''(0)$, and from the geometry of the associated front, that the presence of the acellular gap immediately implies transversal instability of the tumor interface in this model for sufficiently small $\eps>0$. In the benign and malignant no-gap cases, however, the tumor interface can be stable or unstable to long wavelength perturbations depending on parameters. Our characterization of this instability in terms of the coefficient $\lambda_\mathrm{c}''(0)$ allows us to easily determine stability boundaries in parameter space via numerical continuation; see~\S\ref{S:NUM}. We also explore the effect of individual parameters, such as the parameter $\delta_1$, which turns out to be a natural parameter to transition between the three cases of benign and malignant no-gap/gap fronts, and the Allee parameter $a$. In particular, for the latter, we find that in general a stronger Allee effect, i.e. larger $a$, leads to a decrease in the speed of the associated front, but also leads to the onset of the transversal instability.

The remainder of the paper is organized as follows. In~\S\ref{S:E}, we analyze the traveling wave equation associated with~\eqref{eq:GG} to describe the existence construction of traveling fronts in the singular limit $\eps\to0$. Using formal singular perturbation arguments we examine the transversal long wavelength  (in)stability of these fronts in~\S\ref{S:STAB}, and we conclude with numerical simulations and a discussion in~\S\ref{S:NUM}.


\section{Existence of traveling fronts}\label{S:E}


In this section, we show that \eqref{eq:GG} supports traveling fronts. We will largely follow the geometric singular perturbation theory approach of \cite{davis2022traveling}. In \cite{davis2022traveling} the formal asymptotic results of \cite{holder2014model} on the original nondimensionalized Gatenby-Gawlinski model were proven rigorously and a geometric interpretation of the benign and malignant cases, as well as the acellular gap, were given. As the focus of the current manuscript is largely on the stability of the traveling fronts, see \S\ref{S:STAB}, and since the derivation of the existence results for the current setting largely mimics the approaches and proofs of \cite{davis2022traveling}, we only succinctly derive the results (and we refer to \cite{davis2022traveling} for the rigorous constructions in a similar setting).      

We set $(U,V,W)(x,y,\tau) = (u,v,w)(x+\eps^\nu c\tau)$ and search for traveling waves in the traveling coordinate $\xi:= x+\eps^\nu c\tau$. This results in a singularly perturbed traveling wave ordinary differential equation (ODE) 
\begin{align}
\begin{split}\label{eq:GG_twODE}
\eps^\nu cu' &=u(1-u)-\delta_1 uw\,,\\
\eps^\nu cv'  &= \rho v(1-v)(v-a)-\delta_2 vw+((1+\kappa -u) v')'\,,\\
\eps^\nu cw' &= \delta_3(v-w)+\frac{1}{\eps^2}w''\,,
\end{split}
\end{align}
where $'$ means differentiation with respect to $\xi$.
Upon introducing\footnote{See Remark 2.1 of \cite{davis2022traveling} for the rationale behind this non-standard scaling.} $q:= \eps ^{-\nu}(1+\kappa -u) v' - cv$ and $p= w'/\eps$, we rewrite this system as a first order slow/fast ODE
\begin{align}
\begin{split}\label{eq:fast}
\eps^\nu u' &=\frac{u}{c}\left(1-u-\delta_1w\right)\,,\\
v' &= \eps ^\nu \frac{q+c v}{1+\kappa-u}\,,\\
 \eps ^\nu q'  &= -\rho v(1-v)(v-a)+\delta_2 vw\,,\\
w' &= \eps p\,,\\
p' &= \eps\left(c\eps^{\nu+1} p-\delta_3(v-w)\right)\,.
\end{split}
\end{align}
Observe that the $v$-equation of \eqref{eq:fast} is not singular due to the regularisation term $\kappa$ and the fact that we are looking for traveling fronts with normal cell density $u$ between $0$ and $1$. 
The fixed points
\begin{align*}
p_1 = (0,0,0,0,0), \quad p_2 = (1,0,0,0,0), \quad p_3^\pm = (1-\delta_1 V^\pm, V^\pm,c V^\pm, V^\pm,0), \quad p_4^\pm=(0,V^\pm, c V^\pm ,V^\pm,0),
\end{align*}
of~\eqref{eq:fast} in $(u,v,q,w,p)$-space correspond to the steady states $P_1,P_2,P_3^\pm, P_4^\pm$~\eqref{eq:steadystates} of~\eqref{eq:GG}.

We observe three critical $\nu$-cases in~\eqref{eq:fast}: $\nu=0$ and $\nu=\pm 1$. 
In \cite{holder2014model} it was shown for the original nondimensionalized Gatenby-Gawlinski model that, when translated to the current setting, the case $\nu=1$ does not lead to traveling fronts. In contrast, the case $\nu=-1$ led to {\emph{fast traveling fronts}}, while the case $\nu=0$ led to {\emph{slow traveling fronts}}.
In \cite{davis2022traveling} both the fast traveling fronts and slow traveling fronts were investigated. Based on our numerical simulations it appears that the $\nu=0$-case is the most relevant for the instabilities we want to study. Therefore, we focus on the $\nu=0$-case and only briefly highlight the $\nu=-1$-case in the remark below.  

\begin{figure}[ht]
\centering
\includegraphics[width=0.6\linewidth]{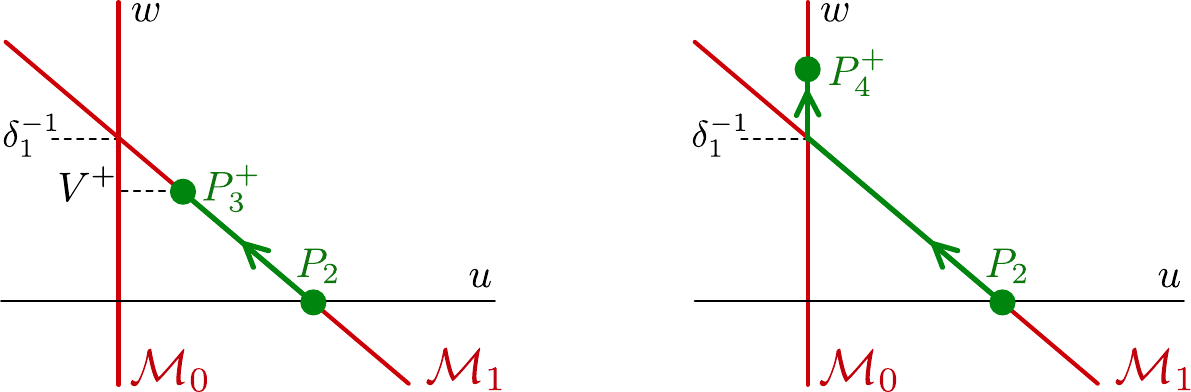}
\caption{Schematic of the singular heteroclinic orbit in the benign case where $u \neq 0$ after the invasion of the traveling wave (left) and the malignant case where $u = 0$ after the invasion (right).}
\label{f:front-schematic}
\end{figure}

For $\nu=0$ in \eqref{eq:fast} we observe that $(u,v,q)$ are {\emph{fast variables}}, while $(w,p)$ are {\emph{slow variables}}.
Letting $\eps \to0$, we find the critical manifolds
\begin{align}
\begin{split}
\label{eq:crit}
\mathcal{M}_0 &= \left\{(u,v,q,w,p) \quad | \quad u=0,\,\, q=-cv,\,\, \rho v(1-v)(v-a) = \delta_2vw  \right\}\,,\\
\mathcal{M}_1 &= \left\{ (u,v,q,w,p) \quad | \quad u=1-\delta_1 w,\,\, q=-cv,\,\, \rho v(1-v)(v-a) = \delta_2 vw  \right\}\,,
\end{split}
\end{align}
which meet along the nonhyperbolic transcritical singularity curve at $w = 1/\delta_1$. The fixed points $p_1,p_4^\pm$ lie on $\mathcal{M}_0$, while $p_2, p_3^\pm$ lie on $\mathcal{M}_1$. In the benign case $\delta_1 V^+<1$, see \eqref{Vplusmin}, we can construct a singular heteroclinic orbit between $p_2$ and $p_3^+$ in the subspace $u=1-\delta_1 w$, that is, the orbit is entirely contained within $\mathcal{M}_1$ and avoids the transcritical singularity; see Fig.~\ref{f:front-schematic}. In the malignant case $\delta_1 V^+>1$, we can construct a singular heteroclinic orbit between $p_2$ and $p_4^+$, but it necessarily passes through the transcritical curve in order to transition from $\mathcal{M}_1$ to $\mathcal{M}_0$. 

In order to construct singular heteroclinic orbits in each case, we consider the associated layer problems describing the dynamics near the interface (the \emph{fast field}) and reduced problems describing the dynamics away from the interface (the \emph{slow fields}) in the next sections.

\begin{Remark}
For the fast traveling fronts with $\nu=-1$, the singularly perturbed problem reduces, upon rescaling space once more, at leading order to  
\begin{align*}
\begin{split}
\dot{u} &=\frac{u}{c}\left(1-u-\delta_1w\right)\,,\\
  \dot{q}  &= - \rho v(1-v)(v-a)+\eps^2 \delta_2 vw\,,\\
\dot{w} &=  p\,,\\
\dot{p} &= c p-\delta_3(v-w)\,,
\end{split}
\end{align*}
where $\dot{}$ denotes differentiation with respect to $\bar{\xi}:=\eps x+ct$, 
on the attracting four dimensional critical manifold 
$$
\mathcal{M}_0^F := \left\{(u,v,q,w,p) \,\, \left| \,\, v = -\dfrac{q}{c}\right.\right\}\,.
$$
That is, the existence of fast traveling fronts boils down to showing the existence of heteroclinic solutions to
\begin{align*}
\begin{split}
\dot{u} &=\frac{u}{c}\left(1-u-\delta_1w\right)\,,\\
  c \dot{v}  &= \rho v(1-v)(v-a)-\delta_2 vw\,,\\
\ddot{w} - c \dot{w} &= -\delta_3(v-w)\,.
\end{split}
\end{align*} 
In this manuscript, we do not pursue this direction any further; however, see Theorem 1.1 of \cite{davis2022traveling}. 
\end{Remark}

\subsection{Layer problem}
We consider the layer problem describing the dynamics near the interface (the fast field) for $\nu=0$. That is, we set $\nu=0$ in \eqref{eq:fast} and take the singular limit $\eps \to 0$ to obtain
\begin{align}
\begin{split}\label{eq:layer}
u' &=\dfrac{u}{c}\left(1-u-\delta_1w\right)\,,\\
v' &= \dfrac{q+cv}{1+\kappa-u}\,,\\
q'  &= -\rho v(1-v)(v-a)+\delta_2 vw\,,\\
\end{split}
\end{align}
where we recall that $'$ denotes differentiation with respect to $\xi=x+ct$,
with the fast variables $w$ and $p$ fixed constants.
The fixed points of \eqref{eq:layer} are related to the critical manifolds $\mathcal{M}_{0,1}$ \eqref{eq:crit}, and note that these have several branches 
\begin{align*}
\mathcal{M}^0_0 &= \left\{(u,v,q,w,p) \,\, | \,\, u=v=q=0\right\}\,,\\
\mathcal{M}^\pm_0 &= \left\{(u,v,q,w,p) \,\, | \,\, u=0,\,\, v = v^\pm(w),\,\, q=-c v^\pm(w) \right\}\,,\\
\mathcal{M}^0_1 &= \left\{(u,v,q,w,p) \,\, | \,\,  u=1-\delta_1 w,\,\, v=q=0\right\}\,,\\
\mathcal{M}^\pm_1 &= \left\{(u,v,q,w,p) \,\, | \,\, u=1-\delta_1 w,\,\, v = v^\pm(w),\,\, q=- c v^\pm(w)  \right\}\,,
\end{align*}
where
\begin{align*}
v^\pm(w) = \frac{1+a\pm\sqrt{(1-a)^2-4\delta_2 w/\rho}}{2}\,.
\end{align*}
We note that the fixed point $p_2$ lies on $\mathcal{M}^0_1$, while $p_3^+$ lies on $\mathcal{M}^+_1$, and $p_4^+$ lies on $\mathcal{M}_0^+$. By inspecting the Jacobian matrix of \eqref{eq:layer}
\begin{align*}
    J = \begin{pmatrix} \frac{1}{c}\left(1-2u-\delta_1 w\right) & 0 & 0\\ \frac{q+cv}{(1+\kappa-u)^2} &\frac{c}{1+\kappa-u} & \frac{1}{1+\kappa-u} \\ 0 & \rho(3v^2-2av+a) +\delta_2 w&0   \end{pmatrix}
\end{align*}
and noting the lower block triangular structure of this matrix, we see that the critical manifolds all lose normal hyperbolicity along the transcritical curve $w=\delta_1^{-1}$ where $1-2u-
\delta_1 w=0$ and along the fold curve where $v^+(w)=v^-(w)$, or equivalently $\rho(3v^2-2av+a) +\delta_2 w=0$. Away from these curves, the manifolds $\mathcal{M}^0_0$, $\mathcal{M}^0_1$, $\mathcal{M}^+_0$, and $\mathcal{M}^+_1$ are all normally hyperbolic and of saddle type in the $(v,q)$-subsystem. The manifolds $\mathcal{M}^-_0$ and $\mathcal{M}^-_1$ are of center type when $c=0$ and normally repelling/attracting in the $(v,q)$-subsystem when $c\lessgtr0$; however these two manifolds will not be important for the analysis.

\begin{figure}[h]
\centering
\includegraphics[width=0.4\linewidth]{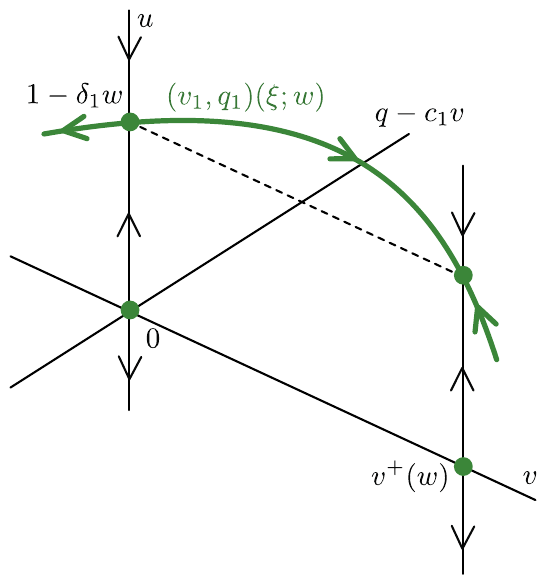}\hspace{0.1\linewidth}
\includegraphics[width=0.4\linewidth]{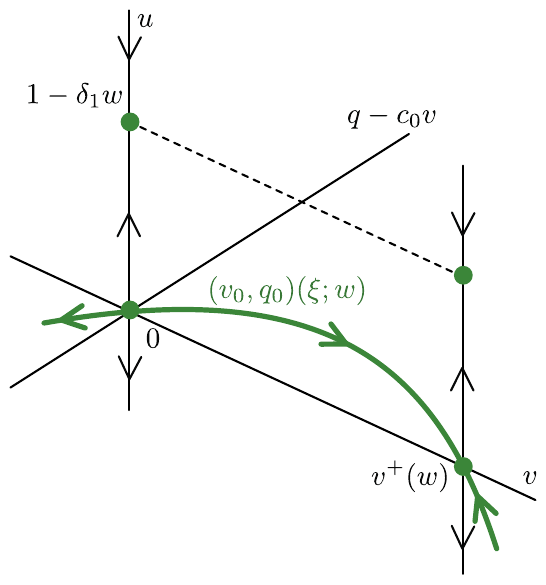}
\caption{(Left) The fast orbit $(v_1,q_1)(\xi;w)$ in the subspace $u=1-\delta_1 w$ with $0 \leq w < \min\{\rho (1-a)^2/(4\delta_2),1/\delta_1\}$. (Right) The fast orbit $(v_0,q_0)(\xi;w)$ in the subspace $u=0$ with $w > 0$.}
\label{f:layer-front}
\end{figure}

Importantly, for the benign case $u=1-\delta_1 w$ in the layer dynamics and
for each $0\leq w < \min\{\rho (1-a)^2/(4\delta_2),1/\delta_1\}$ (within this subspace $u=1-\delta_1 w$) there exists a heteroclinic orbit $(v,q)=(v_1,q_1)(\xi;w)$ between $\mathcal{M}^0_1$ and $\mathcal{M}^+_1$ satisfying the planar ODE
\begin{align*}
\begin{split}
v' &= \frac{q+cv}{\kappa+\delta_1 w}\,,\\
q'  &= -\rho v(1-v)(v-a)+\delta_2 vw\,,\\
\end{split}
\end{align*}
with $c = c_1(w)$, see Fig.~\ref{f:layer-front}. Here, $v_1(\xi;w)$ and $c_1(w)$ are given by
\begin{align*}
    v_1(\xi;w)&:= \frac{v^+(w)}{2}\left(1+\tanh\left(\frac{v^+(w)}{2}\sqrt{\frac{\rho}{2(\kappa+\delta_1w)}} \xi \right)  \right)\,,
\\
c_1(w)&:=\sqrt{\frac{2(\kappa+\delta_1w)}{\rho}}\left(\frac{v^+(w)}{2}-v^-(w)\right).
\end{align*}

Similarly, in the subspace $u=0$, there exists a heteroclinic orbit $(v,q)=(v_0,q_0)(\xi;w)$ between $\mathcal{M}^0_0$ and $\mathcal{M}^+_0$ satisfying the planar ODE
\begin{align*}
\begin{split}
v' &= \frac{q+cv}{1+\kappa}\,,\\
q'  &= -\rho v(1-v)(v-a)+\delta_2 vw\,,\\
\end{split}
\end{align*}
for $c = c_0(w)$, where
\begin{align*}
    v_0(\xi;w)&:= \frac{v^+(w)}{2}\left(1+\tanh\left(\frac{v^+(w)}{2}\sqrt{\frac{\rho}{2(1+\kappa)}} \xi \right)  \right)\,,
\\
c_0(w)&:=\sqrt{\frac{2(1+\kappa)}{\rho}}\left(\frac{v^+(w)}{2}-v^-(w)\right).
\end{align*}

\subsection{Reduced problem}
\begin{figure}[ht]
\centering
\includegraphics[width=0.3\linewidth]{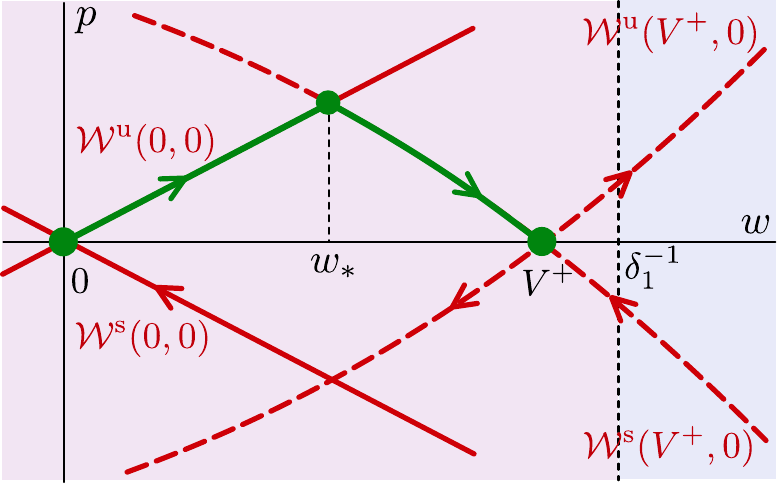}\hspace{0.025\linewidth}
\includegraphics[width=0.3\linewidth]{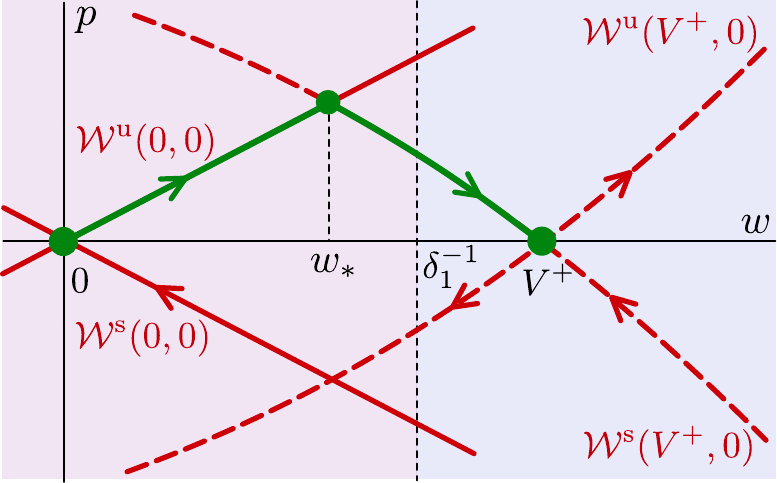}\hspace{0.025\linewidth}
\includegraphics[width=0.3\linewidth]{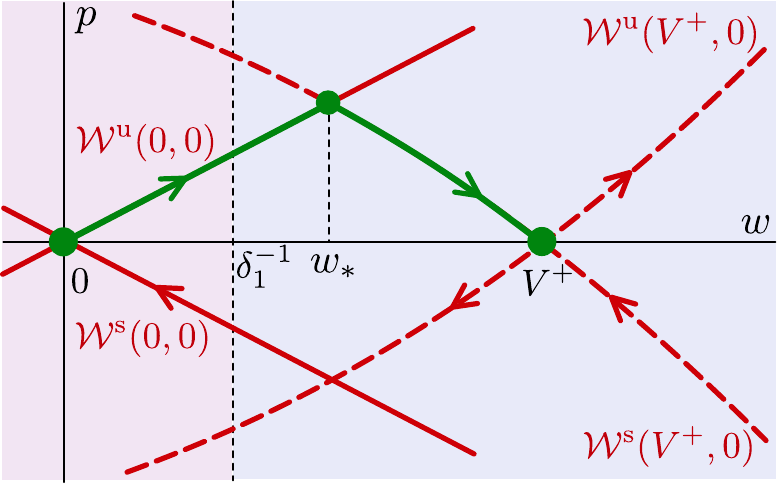}
\caption{The slow orbits on $\mathcal{M}^0_{0,1}$ and $\mathcal{M}^+_{0,1}$, corresponding to $\mathcal{W}^\mathrm{u}(0,0)$ and $\mathcal{W}^\mathrm{u}(V^+,0)$, respectively, in the benign (left), malignant no-gap (center), and malignant gap (right) cases. The manifolds $\mathcal{W}^\mathrm{s,u}(0,0)$ on $\mathcal{M}^0_{0,1}$ are depicted in solid red, while the manifolds $\mathcal{W}^\mathrm{s,u}(V^+,0)$ on $\mathcal{M}^+_{0,1}$ are depicted in dashed red, and the vertical dashed line indicates the subspace $w=\delta_1^{-1}$ at which the reduced dynamics transition from $\mathcal{M}^{0,+}_{0}$ to $\mathcal{M}^{0,+}_{1}$; see also Fig.~\ref{f:3D_singular}. }
\label{f:reduced}
\end{figure}

Rescaling $\zeta=\eps \xi$ and setting $\eps=0$, we obtain the reduced problem describing the dynamics away from the interface (the slow fields)
\begin{align*}
\begin{split}
w_\zeta &= p\,,\\
p_\zeta &= -\delta_3(v-w)\,,
\end{split}
\end{align*}
which does not explicitly depend on $u$ and describes the leading order dynamics on the critical manifolds. When restricted to $\mathcal{M}^0_0$ or $\mathcal{M}_1^0$, where $v=0$, this results in the system
\begin{align}
\begin{split}\label{eq:redM0}
w_\zeta &= p\,,\\
p_\zeta &= \delta_3 w\,,
\end{split}
\end{align}
while on $\mathcal{M}^+_0$ and $\mathcal{M}^+_1$, where $v=v^+(w)$, we have the system
\begin{align}
\begin{split}\label{eq:redM+}
w_\zeta &= p\,,\\
p_\zeta &= \delta_3(w-v^+(w))\,.
\end{split}
\end{align}
The system~\eqref{eq:redM0} admits a saddle equilibrium at $(0,0)$ corresponding to the fixed point $p_2$ of the full system, while~\eqref{eq:redM+} admits a saddle equilibrium at $(V^+,0)$, see \eqref{Vplusmin}, corresponding to $p_3^+$ or $p_4^+$.
The (un)stable manifolds $\mathcal{W}^\mathrm{s,u}(0,0)$ of the equilibrium $(0,0)$ of~\eqref{eq:redM0} within $\mathcal{M}^0_0$ or $\mathcal{M}^0_1$ are given by the lines $\{p = \pm \sqrt{\delta_3}w\}$, while the (un)stable manifolds $\mathcal{W}^\mathrm{s,u}(V^+,0)$ of the equilibrium $(V^+,0)$ of~\eqref{eq:redM+} within $\mathcal{M}^+_0$ or $\mathcal{M}^+_1$ can be determined (implicitly) from the Hamiltonian structure of \eqref{eq:redM+}. In particular, the quantities 
\begin{align*}
    E_0(w,p):=\frac{1}{2}p^2-\frac{\delta_3}{2}w^2\,,
\end{align*}
and
\begin{align*}
    E_+(w,p):=
     E_0(w,p) + E(w) =
    \frac{1}{2}p^2-\frac{\delta_3}{2}w^2+\frac{\delta_3}{2}(V^+)^2+\int_{V^+}^w\delta_3v^+(s)\mathrm{d}s\,,
\end{align*}
are conserved in~\eqref{eq:redM0} and~\eqref{eq:redM+}, respectively, and satisfy
\begin{align*}
    E_0(0,0)=0=E_+(V^+,0)\,.
\end{align*}
The equations \eqref{eq:redM0} and \eqref{eq:redM+} only differ by the term $\delta_3 v^+(w)$ in $p$-component and this term is always strictly positive. Furthermore, the flow of~\eqref{eq:redM+} through the line segment $\{(w,p)=(0,0)+t(V^+,0)\,|\,t \in (0,1)\}$ is downwards, while the flow through $\{(w,p)=(V^+,0)+t( 0,\sqrt{\delta_3} V^+)\,|\,t \in (0,1)\}$ points to the right.
Hence the projection of the unstable manifold $\mathcal{W}^\mathrm{u}(0,0)$ from $\mathcal{M}^0_1$ onto $\mathcal{M}^+_1$ transversely intersects the stable manifold $\mathcal{W}^\mathrm{s}(V^+,0)$ of the equilibrium $(V^+,0)$ at some $(w,p) = (w_*, p_*)$, where $p_* = \sqrt{\delta_3}w_*$, and $0<w_*<V^+$ satisfies
\begin{align*}
    E_+(w_*,\sqrt{\delta_3}w_*)=0\,,
\end{align*}
or equivalently,
\begin{align}\label{eq:wstar_def}
    0=(V^+)^2+\int_{V^+}^{w_*}\left(1+a+\sqrt{(1-a)^2-\frac{4\delta_2 z}{\rho}}\right)\mathrm{d}z \,,
\end{align}
see Fig.~\ref{f:reduced}. The same holds regarding the projection of the unstable manifold $\mathcal{W}^\mathrm{u}(0,0)$ from $\mathcal{M}^0_0$ onto $\mathcal{M}^+_0$, which transversely intersects the stable manifold $\mathcal{W}^\mathrm{s}(W^+,0)$ of the equilibrium $(W^+,0)$ at $w=w_*$.

\subsection{Singular heteroclinic orbits}
Combining orbits from the reduced and layer problems, we can construct singular heteroclinic orbits in the benign ($\delta_1 V^+<1$) and malignant ($\delta_1 V^+>1$) cases. 

\begin{figure}[t]
\centering
\includegraphics[width=0.3\linewidth]{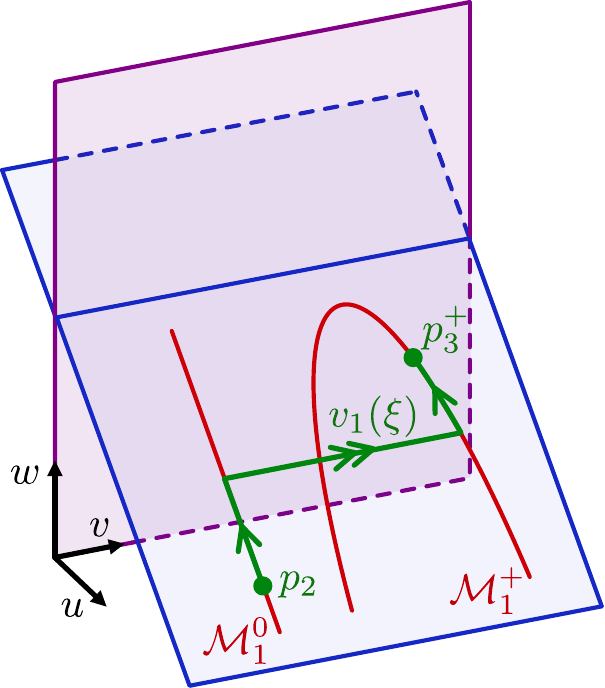}\hspace{0.025\linewidth}
\includegraphics[width=0.3\linewidth]{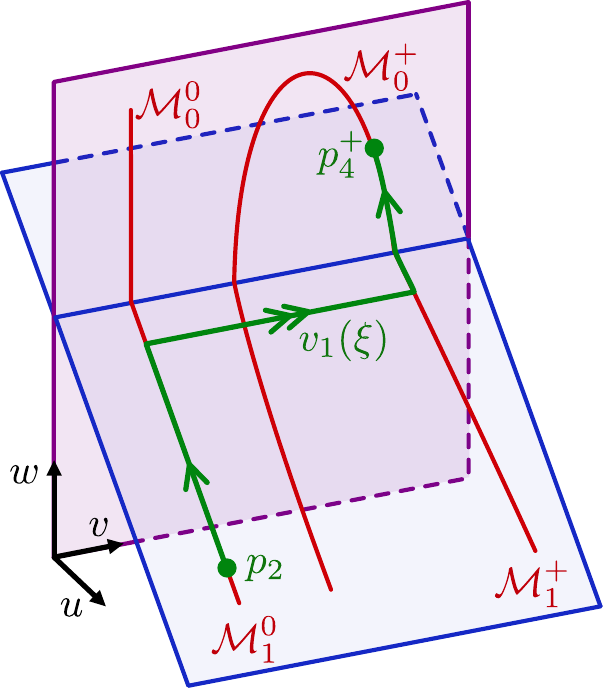}\hspace{0.025\linewidth}
\includegraphics[width=0.3\linewidth]{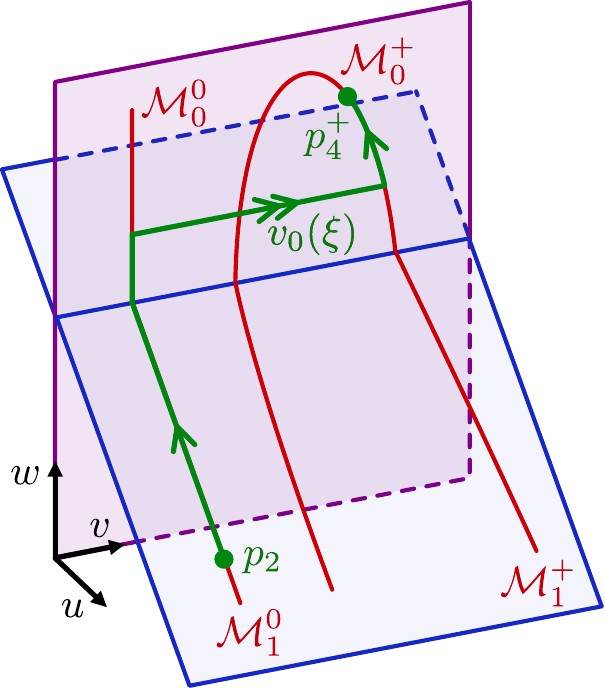}
\caption{Structure of the singular heteroclinic orbit in the benign (left), malignant no-gap (middle), and malignant gap (right) cases. The subspace $\{u=0\}$ is depicted in purple, while the subspace $\{u=1-\delta_1 w\}$ is depicted in blue.}
\label{f:3D_singular}
\end{figure}

\paragraph{Benign case:}In the benign case, we construct a singular heteroclinic orbit between $p_2$ and $p^+_3$ by concatenating the three orbit segments, see Fig.~\ref{f:reduced} (left panel) and Fig.~\ref{f:3D_singular} (left panel):
\begin{enumerate}[1)]
\item slow orbit on $\mathcal{M}^0_1$ given by $\mathcal{W}^\mathrm{u}(0,0)\cap \{0\leq w\leq w_*\}$, with $w^*<V^+$.
\item fast jump $(v_*,q_*)(\xi):=(v_1,q_1)(\xi;w_*)$ with speed $c=c_*(w):=c_1(w_*)$ in the layer problem~\eqref{eq:layer} for $w=w_*$ in the subspace $u = 1-\delta_1w_*$, where $w_*$ is defined as in~\eqref{eq:wstar_def}, see Fig.~\ref{f:layer-front}. Note that $\delta_1w_*<1$ in the benign case.
\item slow orbit on $\mathcal{M}^+_1$ given by $\mathcal{W}^\mathrm{s}(V^+,0)\cap \{w_*< w\leq V^+\}$.
\end{enumerate}
In the benign case, $\mathcal{M}^0_1$ and $\mathcal{M}^+_1$ are normally hyperbolic in the relevant region, and the fast jump is transversely constructed (with the speed $c\approx c_1(w_*)$ as a free parameter). Therefore, we expect that the singular orbit perturbs to a traveling front of the full problem for $0<\eps\ll 1$.

\paragraph{Malignant case:} In the malignant case, we construct a heteroclinic orbit between $p_2$ and $p^+_4$, which necessarily crosses through the transcritical curve $w = \delta_1^{-1}$, see Fig.~\ref{f:front-schematic}. However, we must split this into a further two cases, namely whether the corresponding value of $w_*$ (as defined in~\eqref{eq:wstar_def}) satisfies $w_*\lessgtr \delta_1^{-1}$, as this determines whether the fast jump in the layer problem~\eqref{eq:layer} occurs in the subspace $u=0$ or $u=1-\delta_1w$. That is, whether the fast jump occurs {\emph{before}} of {\emph{after}} crossing the transcritical curve.
In the latter case $w_*>\delta_1^{-1}$, there is a region in which the values of both $u,v$ are zero, \emph{before} the fast jump occurs, that is, the front admits a gap, called the acellular gap, where only acid is present.
In the former case $w_*<\delta_1^{-1}$ there is no acellular gap (as in the benign case above). The boundary between these cases in parameter space is determined by the condition
\begin{align*}
    0=(V^+)^2+\int_{V^+}^{\delta_1^{-1}}1+a+\sqrt{(1-a)^2-\frac{4\delta_2 z}{\rho}}\mathrm{d}z \,,
\end{align*}
where $V^+$ is as in~\eqref{Vplusmin}. We refer to these (sub)cases as the malignant gap and malignant no-gap cases, respectively.

\begin{figure}
\centering
\includegraphics[width=0.4\linewidth]{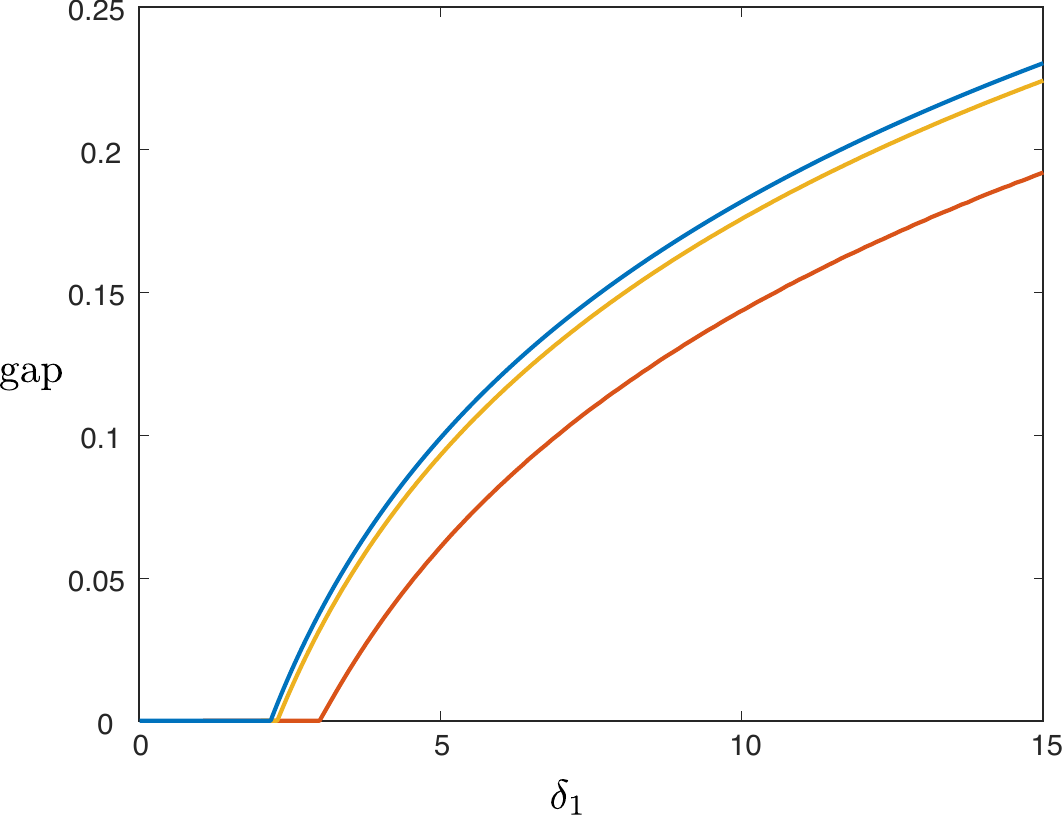}
\caption{Plot of gap width versus $\delta_1 \in(0.1,15)$ obtained by numerical continuation of the traveling wave equation~\eqref{eq:fast} for the parameter values $(a,\kappa, \delta_2,\delta_3, \rho) = (0.1,0.1,0.1,70,1.0)$ and $\eps=0.0063$ (red), $\eps=10^{-4}$ (yellow). The gap width was computed by measuring the spatial width where both the $u$ and $v$ profiles of the corresponding front solution were below a threshold value of $10\eps$. Also plotted (blue) is the singular limit gap width obtained by solving~\eqref{eq:wstar_def} for $w_*$ using Mathematica and integrating~\eqref{eq:redM0} to obtain the time spent along $\mathcal{M}^0_0$ between $w=\delta_1^{-1}$ and $w=w_*$. }
\label{f:gapwidth}
\end{figure}

In the {\bf malignant no-gap case} $w_*<\delta_1^{-1}$, we have the following singular concatenation, see Fig.~\ref{f:reduced} (middle panel) and Fig.~\ref{f:3D_singular} (middle panel):
\begin{enumerate}[1)]
\item slow orbit on $\mathcal{M}^0_1$ given by $\mathcal{W}^\mathrm{u}(0,0)\cap \{0\leq w\leq w_*\}$.
\item fast jump $(v_*,q_*)(\xi):=(v_1,q_1)(\xi;w_*)$ with speed $c=c_*(w):=c_1(w_*)$ in the layer problem~\eqref{eq:layer} for $w=w_*$ in the subspace $u = 1-\delta_1w_*$, where $w_*$ is defined as in~\eqref{eq:wstar_def}, see Fig.~\ref{f:layer-front}. Note that $\delta_1w_*<1$ in the malignant no-gap case.
\item slow orbit on $\mathcal{M}^+_1$ given by $\mathcal{W}^\mathrm{s}(V^+,0)\cap \{w_*< w\leq \delta_1^{-1}\}$.
\item slow orbit on $\mathcal{M}^+_0$ given by $\mathcal{W}^\mathrm{s}(V^+,0)\cap \{\delta_1^{-1}\leq w< V^+\}$.
\end{enumerate}
In the {\bf malignant gap case} $w_*>\delta_1^{-1}$, we have the concatenation, see Fig.~\ref{f:reduced} (right panel) and Fig.~\ref{f:3D_singular} (right panel):
\begin{enumerate}[1)]
\item slow orbit on $\mathcal{M}^0_1$ given by $\mathcal{W}^\mathrm{u}(0,0)\cap \{0\leq w\leq \delta_1^{-1}\}$.
\item slow orbit on $\mathcal{M}^0_0$ given by $\mathcal{W}^\mathrm{u}(0,0)\cap \{\delta_1^{-1}< w\leq w_*\}$.
\item fast jump $(v_*,q_*)(\xi):=(v_0,q_0)(\xi;w_*)$ with speed $c=c_*(w):=c_0(w_*)$ in the layer problem~\eqref{eq:layer} for $w=w_*$ in the subspace $u = 0$, where $w_*$ is defined as in~\eqref{eq:wstar_def}, see Fig.~\ref{f:layer-front}. 
Note that $\delta_1w_*>1$ in this case.
\item slow orbit on $\mathcal{M}^+_0$ given by $\mathcal{W}^\mathrm{s}(V^+,0)\cap \{w_*< w\leq W^+\}$.
\end{enumerate}
The acellular gap manifests as the slow orbit portion  $\mathcal{W}^\mathrm{u}(0,0)\cap \{\delta_1^{-1}\leq w\leq w_*\}$ on $\mathcal{M}^0_0$. In general, one expects the gap size to thus increase with $\delta_1$; see Fig.~\ref{f:gapwidth}.

In both cases (malignant gap and malignant no-gap), the singular orbit traverses the transcritical curve $\delta_1w=1$ along a slow orbit. In the critical crossover case $w_*=\delta_1^{-1}$, the fast jump occurs precisely along the transcritical singularity curve. Due to the lack of hyperbolicity which occurs at the transcritical bifurcation, the persistence of orbits in the malignant case for $0<\eps\ll1$ is nontrivial, necessitating the use of blow-up desingularization methods~\cite{krupa2001extending}. Nevertheless we assume their persistence for sufficiently small $\eps>0$.

\section{Stability of planar tumor interfaces}\label{S:STAB}
Given a traveling front $(u_\mathrm{h},v_\mathrm{h},w_\mathrm{h})(\xi)$ with speed $c_\mathrm{h} = c_*(w_*)+\mathcal{O}(\eps)$ constructed using the slow/fast structure as in the preceding section, we now consider its spectral stability in two space dimensions, focusing on the long wavelength (in)stability criterion as explored in~\cite{CDLOR}. Thus, like in~\cite{CDLOR}, we do not explicitly analyze the stability of the front $(u_\mathrm{h},v_\mathrm{h},w_\mathrm{h})(\xi)$ with respect to longitudinal perturbations, i.e. perturbations that only depend on $\xi$ (or $x$). If the front is longitudinally stable -- as is strongly suggested by our numerical results --  then the upcoming stability criterion determines the (in)stability of the front to long wavelength perturbations in the $y$-direction, i.e. transverse to the direction of propagation. (We again refer to~\cite{CDLOR} for a review of the literature on methods by which the longitudinal (in)stability of a singular traveling front as $(u_\mathrm{h},v_\mathrm{h},w_\mathrm{h})(\xi)$ can be established analytically.) We also refer ahead to~\S\ref{S:NUM} for numerical evidence (see Fig.~\ref{f:profiles-spectra}) that the traveling fronts under consideration are $1$D spectrally stable.

\subsection{Long wavelength (in)stability}
In a comoving frame, we rewrite~\eqref{eq:GG} as 
\begin{align}
\begin{split}\label{eq:GGsco}
U_\tau &=F(U,W)-c U_\xi \,,\\
V_\tau  &= G(U,V,W)+\nabla\cdot((1+\kappa-U)\nabla V)-c V_\xi\,,\\
W_\tau &= H(V,W)+\frac{1}{\eps^2}\Delta W -c W_\xi \,,
\end{split}
\end{align}
where
\begin{align}\label{eq:FGHdef}
F(U,W)&=U(1-U)-\delta_1 UW \,,\\
G(U,V,W)&=\rho V(1-V)(V-a)-\delta_2 VW\,,\\
H(V,W)&=\delta_3(V-W)\,.
\end{align}
and $(U,V,W)=(U,V,W)(\xi,y,\tau)$. Upon substituting the ansatz $(U,V,W)=(u_\mathrm{h},v_\mathrm{h},w_\mathrm{h})(\xi)+(\bar{u},\bar{v},\bar{w})(\xi)e^{i\ell y +\lambda \tau}$ and $c=c_\mathrm{h}$, this results in the linear eigenvalue problem
\begin{align}
\begin{split}\label{eq:GGseval}
\lambda u &=F_{u}u+F_w w-c_\mathrm{h} u_\xi\,,\\
\lambda v  &= G_{v}v+G_w w+(1+k-u_\mathrm{h}){v}_{\xi \xi}-c_\mathrm{h} {v}_\xi-uv_{\mathrm{h},\xi\xi}-{u}_{\xi}v_{\mathrm{h},\xi}-u_{\mathrm{h},\xi}{v}_{\xi}-\ell^2(1+k-u_\mathrm{h}){v}\,,\\
\lambda w &= H_{v}v+H_w w+\frac{1}{\eps^2}w_{\xi\xi}-c_\mathrm{h} w_\xi-\frac{\ell^2}{\eps^2}w \,,
\end{split}
\end{align}
where we have dropped the bars for convenience and we denote $F_{u}(\xi):=\partial_{u}F(u_\mathrm{h}(\xi), w_\mathrm{h}(\xi))$, etc. We write~\eqref{eq:GGseval} in the form
\begin{align}\label{eq:GGseval_mat}
  \mathbb{L} (\xi) \begin{pmatrix}  u \\ v \\ w \end{pmatrix} 
  &= \lambda
  \begin{pmatrix} u \\ v \\ w \end{pmatrix} 
  + \ell^2 
  \begin{pmatrix} 0 \\ (1+\kappa-u_\mathrm{h})v \\ \frac{1}{\eps^2} w \end{pmatrix}\,,
\end{align}
where 
\begin{align*}
\mathbb{L}(\xi):= \begin{pmatrix}-c_\mathrm{h}\partial_\xi+ F_{u} & 0 & F_w \\ -v_{\mathrm{h},\xi\xi}-v_{\mathrm{h},\xi}\partial_\xi & (1+\kappa-u_\mathrm{h})\partial_{\xi\xi}-c_\mathrm{h}\partial_\xi + G_{v}  -u_{\mathrm{h},\xi} \partial_\xi & G_w \\ 0 & H_{v} & \frac{1}{\eps^2}\partial_{\xi\xi}-c_\mathrm{h}\partial_\xi +H_w   \end{pmatrix}\,.
\end{align*}
Due to translation invariance, the derivative $(u,v,w)(\xi) = (u_\mathrm{h}, v_\mathrm{h}, w_\mathrm{h})'(\xi)$ of the wave with respect to $\xi$ satisfies~\eqref{eq:GGseval_mat} when $\lambda=\ell=0$. To examine long-wavelength interfacial instabilities in the direction transverse to the front propagation, we expand about this solution for small $|\ell|\ll1$,  and noting the symmetry $\ell \to -\ell$ we obtain
\begin{align*}
  \lambda_\mathrm{c} (\ell) &= \lambda_{\mathrm{c},2} \ell^2 + \mathcal{O}(\ell^4)\,, \qquad 
  \begin{pmatrix} u(\xi; \ell^2)  \\
  v(\xi; \ell^2)  \\ w(\xi; \ell^2)   \end{pmatrix} 
  =
  \begin{pmatrix}
      u_{\mathrm{h},\xi} (\xi) \\
      v_{\mathrm{h},\xi} (\xi) \\
      w_{\mathrm{h},\xi} (\xi)
  \end{pmatrix}+
  \begin{pmatrix}
      \tilde{u}(\xi) \\
      \tilde{v}(\xi) \\
      \tilde{w}(\xi)
  \end{pmatrix} \ell^2 + \mathcal{O}(\ell^4)\,,
\end{align*}
and substitute into~\eqref{eq:GGseval_mat} to obtain
\begin{align*}
  \mathbb{L} \begin{pmatrix}
      \tilde{u}(\xi) \\
      \tilde{v}(\xi) \\
      \tilde{w}(\xi)
  \end{pmatrix} 
  &= \lambda_{\mathrm{c},2} \begin{pmatrix}
        u_{\mathrm{h},\xi} (\xi) \\
      v_{\mathrm{h},\xi} (\xi) \\
      w_{\mathrm{h},\xi} (\xi)
  \end{pmatrix} +
  \begin{pmatrix}
     0 \\
      (1+\kappa-u_\mathrm{h})v_{\mathrm{h},\xi} (\xi) \\
      \frac{1}{\eps^2}w_{\mathrm{h},\xi} (\xi)
  \end{pmatrix}.
\end{align*}
This results in the Fredholm solvability condition
\begin{align}\label{eq:solvability}
0= \Bigg\langle \lambda_{\mathrm{c},2}\begin{pmatrix}
         u_{\mathrm{h},\xi} (\xi) \\
      v_{\mathrm{h},\xi} (\xi) \\
      w_{\mathrm{h},\xi} (\xi)
  \end{pmatrix}+
  \begin{pmatrix}
     0 \\
      (1+\kappa-u_\mathrm{h})v_{\mathrm{h},\xi} (\xi) \\
      \frac{1}{\eps^2}w_{\mathrm{h},\xi} (\xi)
  \end{pmatrix}, \begin{pmatrix}
      u^A (\xi) \\
      v^A (\xi) \\
      w^A (\xi)
  \end{pmatrix} \Bigg\rangle_{L^2}\,,
\end{align}
where $(u^A,v^A, w^A)(\xi)$ denotes the unique bounded solution to the adjoint equation 
\begin{align}\label{eq:adjointequation}
    \mathbb{L}^A(\xi)\begin{pmatrix}
      u (\xi) \\
     v (\xi) \\
      w (\xi)
  \end{pmatrix}=0\,,
\end{align}
where the adjoint operator $\mathbb{L}^A$ is given by
\begin{align*}
  \mathbb{L}^A (\xi) &=
  \begin{pmatrix}
   c_\mathrm{h} \partial_\xi + F_{u} & v_{\mathrm{h},\xi} \partial_\xi & 0 \\
    0 & (1+\kappa - u_\mathrm{h}) \partial_{\xi\xi} + c_\mathrm{h} \partial_\xi + G_{v} - u_{\mathrm{h},\xi} \partial_\xi & H_{v} \\
    F_w & G_w
    & \frac{1}{\varepsilon^2} \partial_{\xi\xi}+c_\mathrm{h}\partial_\xi + H_w
  \end{pmatrix} .
\end{align*}
Solving~\eqref{eq:solvability} for $\lambda_{\mathrm{c},2}$, we find
\begin{align}\label{eq:lambda2c_expression}
  \lambda_{c,2} &= - \dfrac{\displaystyle\int_\mathbb{R} \left(1+\kappa - u_\mathrm{h}(\xi)\right)v_{\mathrm{h},\xi}(\xi) v^A(\xi) + \frac{1}{\eps^2} w_{\mathrm{h},\xi}(\xi) w^A(\xi) \mathrm{d}\xi }{\displaystyle \int_\mathbb{R} u_{\mathrm{h},\xi}(\xi) u^A(\xi) + v_{\mathrm{h},\xi}(\xi) v^A(\xi) + w_{\mathrm{h},\xi}(\xi) w^A(\xi) \mathrm{d} \xi}.
\end{align}
We note that the sign of $\lambda_{\mathrm{c},2}$ determines the stability of the interface to transverse long wavelength perturbations. To estimate the expression~\eqref{eq:lambda2c_expression}, we need to obtain leading-order approximations of the adjoint solution $(u^A,v^A, w^A)(\xi)$.

\subsection{Leading order asymptotics of $\lambda_{\mathrm{c},2}$}\label{sec:lambda2c_asymp}

We consider the fast formulation
  \begin{align}
  \begin{split}
 c_\mathrm{h}u_\xi +F_{u}(u_\mathrm{h}(\xi),w_\mathrm{h}(\xi))u + v_{\mathrm{h},\xi} v_\xi &=0\,,\\
      (1+\kappa - u_\mathrm{h}) v_{\xi\xi}+ c_\mathrm{h} v _\xi  +G_v(u_\mathrm{h}(\xi),v_\mathrm{h}(\xi),w_\mathrm{h}(\xi)) v+H_v(v_\mathrm{h}(\xi),w_\mathrm{h}(\xi))w - u_{\mathrm{h},\xi} v_\xi 
      &= 0\,, \\
      w_{\xi\xi}+ \eps^2c_\mathrm{h} w_\xi+\varepsilon^2\left( F_w(u_\mathrm{h}(\xi),w_\mathrm{h}(\xi))u + G_w(u_\mathrm{h}(\xi),v_\mathrm{h}(\xi),w_\mathrm{h}(\xi))v + H_w(v_\mathrm{h}(\xi),w_\mathrm{h}(\xi))w  \right) &= 0\,,\\
\end{split}
\label{eq:fastpart}
  \end{align}
of the adjoint equation~\eqref{eq:adjointequation}, and the associated 
slow formulation
    \begin{align}
  \begin{split}
 \eps c_\mathrm{h}u_\zeta +F_{u}(u_\mathrm{h}(\zeta/\eps),w_\mathrm{h}(\zeta/\eps))u + \eps^2 v_{\mathrm{h},\zeta} v_\zeta &=0\,,\\
      \eps^2 (1+\kappa - u_\mathrm{h}) v_{\zeta\zeta}+ \eps c_\mathrm{h} v_\zeta  +G_v(u_\mathrm{h}(\zeta/\eps),v_\mathrm{h}(\zeta/\eps),w_\mathrm{h}(\zeta/\eps)) v+H_v(v_\mathrm{h}(\zeta/\eps),w_\mathrm{h}(\zeta/\eps))w - \eps^2 u_{\mathrm{h},\zeta} v_\zeta 
      &= 0 \,,\\
      w_{\zeta\zeta}+ \eps c_\mathrm{h} w_\zeta+F_w(u_\mathrm{h}(\zeta/\eps),w_\mathrm{h}(\zeta/\eps))u + G_w(u_\mathrm{h}(\zeta/\eps),v_\mathrm{h}(\zeta/\eps),w_\mathrm{h}(\zeta/\eps))v + H_w(v_\mathrm{h}(\zeta/\eps),w_\mathrm{h}(\zeta/\eps))w  &= 0\,,
\end{split}
\label{eq:slowpart}
  \end{align}
   where $\zeta =\eps \xi$, and we abuse notation by writing $u=u(\xi)$ in~\eqref{eq:fastpart}, and $u=u(\zeta)$ in~\eqref{eq:slowpart}, etc. In the fast field near the interface, to leading order we have $c_\mathrm{h}=c_*(w_*)$, $w_\mathrm{h}=w_*$, $v_\mathrm{h}=v_*(\xi)$, and $u_\mathrm{h}=u_*$, where
  \begin{align*}
      c_*(w_*)&=\begin{cases}c_1(w_*),& w_*<\delta_1^{-1}\,,\\ c_0(w_*), &w_*>\delta_1^{-1}\,,\end{cases}\\
      (v_*,q_*)(\xi;w_*)&=\begin{cases}(v_1,q_1)(\xi;w_*), & w_*<\delta_1^{-1}\,,\\ (v_0,q_0)(\xi;w_*), & w_*>\delta_1^{-1}\,,\end{cases}\\
            u_*(w_*)&=\begin{cases}1-\delta_1w_*, & w_*<\delta_1^{-1}\,,\\ 0, &w_*>\delta_1^{-1}\,.\end{cases}
  \end{align*}

Thus, to leading order~\eqref{eq:fastpart} becomes
  \begin{align*}
  \begin{split}
 c_*u_\xi +F_{u}(u_*,w_*)u + v_{*,\xi} v_\xi &=0\,,\\
      (1+\kappa - u_*) v_{\xi\xi}+ c_* v _\xi  +G_v(u_*,v_*(\xi),w_*) v+H_v(v_*(\xi),w_*)w
      &= 0\,, \\
      w_{\xi\xi} &= 0.
\end{split}
  \end{align*}
Hence $w$ is constant to leading order with $w=\bar{w}_*$, and $v$ satisfies
  \begin{align}\label{eq:reduced_adjoint}
      (1+\kappa - u_*) v_{\xi\xi}+ c_* v _\xi  +G_v(u_*,v_*(\xi),w_*) v &= -H_v(v_*(\xi),w_*)\bar{w}_*\,.
  \end{align}
The unique bounded solution of the fast reduced adjoint equation
  \begin{align*}
     \mathcal{L}_v^Av= (1+\kappa - u_*) v_{\xi\xi}+ c_* v _\xi  +G_v(u_*,v_*(\xi),w_*) v &= 0
  \end{align*}
  is given by 
   $\bar{v}_*(\xi):=v_{*,\xi}(\xi)e^{-c_* \xi /(1+\kappa-u_*)}$. 
  Since $v_{*,\xi}$ lies in the kernel of $\mathcal{L}_v$,~\eqref{eq:reduced_adjoint} implies that
    \begin{align*}
     0= \bar{w}_*\int_\mathbb{R} H_v(v_*(\xi),w_*)v_{*,\xi}(\xi)\mathrm{d}\xi= \bar{w}_*\left(H(v_*^+,w_*)-H(0,w_*)\right)= \delta_3\bar{w}_*v_*^+,
  \end{align*}
  where $v_*^+:= v^+(w_*)= \lim_{\xi\to\infty}v_*(\xi)$ and we recall that $\lim_{\xi \to -\infty} v_*(\xi)=0$. Hence to leading order $\bar{w}_*=0$ and $v(\xi)=\alpha_*\bar{v}_*(\xi)$. From this we see that $u$ satisfies 
    \begin{align*}
c_*u_\xi +F_{u}(u_*,w_*)u + \alpha_*v_{*,\xi} \bar{v}_{*,\xi}(\xi) &=0
  \end{align*}
to leading order. Since $F_{u}(u_*,w_*) = 1-2u_*-\delta_1w_* < 0$ (both for $\delta_1 w_* \lessgtr 1$) and since $u(\xi)$ must be bounded, it follows that
    \begin{align*}
u(\xi) &=\alpha_*\bar{u}_*(\xi):=-\frac{\alpha_*}{c_*}\int_{\infty}^\xi e^{-\frac{1}{c_*}(1-2u_*-\delta_1w_*) (\xi-s)} v_{*,s}(s) \bar{v}_{*,s}(s)\mathrm{d}s.
  \end{align*}

In the slow fields away from the interface, to leading order 
\begin{align}
    \left(u_\mathrm{h}(\zeta/\eps), v_\mathrm{h}(\zeta/\eps), w_\mathrm{h}(\zeta/\eps)\right)= \left( u_*(w^\pm(\zeta)),f^\pm(w^\pm(\zeta)), w^\pm(\zeta)\right),
\end{align} where $w^-(\zeta)$ denotes the slow orbit of~\eqref{eq:redM0} on $\mathcal{M}^0_0\cup\mathcal{M}^0_1$ corresponding to $\mathcal{W}^\mathrm{u}(0,0)$ and satisfies $w^-(0)=w_*$, and similarly $w^+(\zeta)$ denotes the slow orbit of~\eqref{eq:redM+} on $\mathcal{M}^+_0\cup\mathcal{M}^+_1$ corresponding to $\mathcal{W}^\mathrm{s}(W^+,0)$ and satisfies $w^+(0)=w_*$. Note that whether $w^-$ is contained entirely within $\mathcal{M}^0_1$, or whether $w^+$ is contained entirely within $\mathcal{M}^+_1$, depends on which case we are in, namely, the benign, malignant gap or no-gap cases. To simplify the following computations, we use the analogous notation $f^\pm(w)$ to define the corresponding $v$ coordinate along the slow orbits by
\begin{align*}
    f^-(w) &:= 0, \qquad f^+(w) :=v^+(w),
\end{align*}
and we recall that 
\begin{align*}
               u_*(w)&=\begin{cases}1-\delta_1w, & w<\delta_1^{-1}\\ 0, &w>\delta_1^{-1}\end{cases}.
\end{align*}

Thus, to leading order~\eqref{eq:slowpart} becomes
\begin{align*}
  \begin{split}
F_{u}(u_*(w^\pm),w^\pm)u &=0\,,\\
      G_v(u_*(w^\pm),f^\pm(w^\pm),w^\pm) v+H_v(f^\pm(w^\pm),w^\pm)w &= 0\,, \\
      w_{\zeta\zeta}+F_w(u_*(w^\pm),w^\pm)u + G_w(u_*(w^\pm),f^\pm(w^\pm),w^\pm)v + H_w(f^\pm(w^\pm),w^\pm)w  &= 0\,,
\end{split}
  \end{align*}
from which we deduce that $u=0$, and 
\begin{align*}
      v&=-\frac{H_v(f^\pm(w^\pm),w^\pm)}{G_v(u_*(w^\pm),f^\pm(w^\pm),w^\pm) }w\,,
  \end{align*}
so that $w$ satisfies
\begin{align}\label{eq:reduced_adjoint_slow}
      w_{\zeta\zeta} + \left(H_w(f^\pm(w^\pm),w^\pm)-\frac{H_v(f^\pm(w^\pm),w^\pm)}{G_v(u_*(w^\pm),f^\pm(w^\pm),w^\pm) }G_w(u_*(w^\pm),f^\pm(w^\pm),w^\pm) \right)w &= 0.
\end{align}
Noting that 
\begin{align*}
    (f^\pm)'(w^\pm)=-\frac{G_w(u_*(w^\pm),f^\pm(w^\pm),w^\pm)}{G_v(u_*(w^\pm),f^\pm(w^\pm),w^\pm) }\,,
\end{align*}
the equation~\eqref{eq:reduced_adjoint_slow} becomes
\begin{align*}
      w_{\zeta\zeta} + \left(H_w(f^\pm(w^\pm),w^\pm)+H_v(f^\pm(w^\pm),w^\pm)(f^\pm)'(w^\pm) \right)w &= 0,
\end{align*}
from which we deduce that $w(\zeta) = \alpha^\pm w^\pm_\zeta(\zeta)$ in the slow fields.

We recall that in the fast field $\bar{w}_*=0$ to leading order so that $w=\mathcal{O}(\eps)$, and hence in the slow fields we take $w = \eps \bar{\alpha}^\pm w^\pm_\zeta=\bar{\alpha}^\pm w^\pm_\xi$. Since $w^+_\zeta(0)=w^-_\zeta(0)=p_*$, to ensure continuity, we take $\bar{\alpha}^+=\bar{\alpha}^-=\bar{\alpha}$, so that $\bar{w}_*=\eps\bar{\alpha}p_*$. The jump in $w_\xi$ in the slow fields
\begin{align*}
  \Delta_s w_\xi &= \lim_{\zeta \downarrow 0} w_\xi - \lim_{\zeta \uparrow 0} w_\xi = \bar{\alpha} \left[ \lim_{\zeta \downarrow 0} w^+_{\xi\xi} - \lim_{\zeta \uparrow 0 } w^-_{\xi\xi} \right] = \varepsilon^2 \bar{\alpha} \left[ \lim_{\zeta \downarrow 0} w_{\zeta\zeta}^+ - \lim_{\zeta \uparrow 0} w_{\zeta\zeta}^- \right] \\
  &= - \varepsilon^2 \bar{\alpha} \left[ H(v_*^+, w_*) -
  H( 0, w_*) \right] = - \varepsilon^2 \bar{\alpha} \delta_3 v_*^+
\end{align*}
must be accounted for by the change over the fast field
\begin{align*}
    \Delta_f w_\xi &= -\varepsilon^2 \alpha_* \int_\mathbb{R} F_w(u_*,w_*)\bar{u}_*(\xi) + G_w (u_*,v_*(\xi), w_*) \bar{v}_*(\xi) \mathrm{d} \xi\,,
\end{align*}
from which we obtain
\begin{align}
  \frac{\bar{\alpha}}{\alpha_*} &=
  \frac{\int_\mathbb{R} F_w(u_*,w_*)\bar{u}_*(\xi) + G_w (u_*,v_*(\xi), w_*) \bar{v}_*(\xi) \mathrm{d} \xi}{\delta_3 v_*^+}.
\label{eq:alpharatio}
\end{align}

We can now estimate~\eqref{eq:lambda2c_expression}. At leading order
\begin{align*}
  \int_\mathbb{R} w_{\mathrm{h},\xi}(\xi) w^A(\xi) \mathrm{d} \xi &=
  \left( \int_{-\infty}^{-\frac{1}{\sqrt{\varepsilon} }} + \int_{-\frac{1}{\sqrt{\varepsilon}} }^{\frac{1}{\sqrt{\varepsilon}}} + \int_{\frac{1}{\sqrt{\varepsilon} } }^\infty \right) w_{\mathrm{h},\xi}^*(\xi) w^A(\xi) \mathrm{d} \xi
\\
&= \overline{\alpha}
\int_{-\infty}^{-\frac{1}{\sqrt{\varepsilon} }} (w_\xi^-)^2 \mathrm{d} \xi +
  \varepsilon^2 \bar{\alpha}
  \int_{-\frac{1}{\sqrt{\varepsilon} }}^{\frac{1}{\sqrt{\varepsilon} }} (q_*)^2 \mathrm{d} \xi
+\overline{\alpha}
\int_{\frac{1}{\sqrt{\varepsilon} }}^\infty (w_\xi^+)^2 \mathrm{d} \xi \\
&=  \varepsilon \bar{\alpha} \left( \int_{-\infty}^0 (w_\zeta^-)^2 \mathrm{d} \zeta + \int_0^\infty (w_\zeta^+)^2 \mathrm{d} \zeta \right) + \mathcal{O}(\varepsilon \sqrt{\varepsilon} )\,,
\end{align*}
and
\begin{align*}
  \int_\mathbb{R} v_{\mathrm{h},\xi}(\xi) v^A(\xi) \mathrm{d} \xi
  &= \alpha_* \int_\mathbb{R}  v_{*,\xi}(\xi)^2 
  e^{-c_* \xi/(1+\kappa - u_*)} 
  \mathrm{d} \xi + \mathcal{O}(\varepsilon).
\end{align*}
We note that $u_{\mathrm{h},\xi}=\mathcal{O}(\eps)$ in the fast field while $u^A(\xi)=0$ to leading order in the slow fields, so by~\eqref{eq:lambda2c_expression} and \eqref{eq:alpharatio}, we have that
\begin{align}\label{eq:lambda2cfinal}
  \lambda_{c,2} &\sim 
  - \frac{\bar{\alpha}}{\eps \alpha_*}   \frac{ \displaystyle \left( \int_{-\infty}^0 (w_\zeta^-)^2 \mathrm{d} \zeta + \int_0^\infty (w_\zeta^+)^2 \mathrm{d} \zeta \right) }{\displaystyle\int_\mathbb{R} v_{*,\xi}(\xi)^2 e^{-\left( \frac{c_*}{1+\kappa - u_*} \right)\xi}\mathrm{d} \xi }\nonumber \\
  &=- \frac{1}{\eps \delta_3 v_*^+}\int_\mathbb{R} F_w(u_*,w_*)\bar{u}_*(\xi) + G_w (u_*,v_*(\xi), w_*) \bar{v}_*(\xi) \mathrm{d} \xi \frac{\displaystyle \left( \int_{-\infty}^0 (w_\zeta^-)^2 \mathrm{d} \zeta + \int_0^\infty (w_\zeta^+)^2 \mathrm{d} \zeta \right) }{\displaystyle \int_\mathbb{R} v_{*,\xi}(\xi)^2 e^{-\left( \frac{c_*}{1+\kappa - u_*} \right)\xi}\mathrm{d} \xi }\nonumber \\
  &= \frac{1}{\eps \delta_3 v_*^+}\frac{\displaystyle \left( \int_{-\infty}^0 (w_\zeta^-)^2 \mathrm{d} \zeta + \int_0^\infty (w_\zeta^+)^2 \mathrm{d} \zeta \right) }{\displaystyle \int_\mathbb{R} v_{*,\xi}(\xi)^2 e^{-\left( \frac{c_*}{1+\kappa - u_*} \right)\xi}\mathrm{d} \xi } \left( \int_\mathbb{R} \delta_1u_*\bar{u}_*(\xi) +\delta_2v_*(\xi) \bar{v}_*(\xi) \mathrm{d} \xi\right)
\end{align}
at leading order in $\varepsilon$, where
    \begin{align*}
    \bar{v}_*(\xi)&=v_{*,\xi}(\xi)e^{-\frac{c_*}{1+\kappa-u_*}\xi}\\
\bar{u}_*(\xi)&=\frac{1}{c_*}\int^{\infty}_\xi e^{-\frac{1}{c_*}(1-2u_*-\delta_1w_*) (\xi-s)} v_{*,s}(s) \bar{v}_{*,s}(s)\mathrm{d}s.
  \end{align*}
In particular, the sign of $\lambda_{\mathrm{c},2}$ is determined at leading order by
\begin{align}\label{eq:lambda2c_sgn}
  \mathrm{sign}\left(\lambda_{c,2}\right)&= \mathrm{sign}\left(\int_\mathbb{R} \delta_1u_*\bar{u}_*(\xi) +\delta_2v_*(\xi) \bar{v}_*(\xi) \mathrm{d} \xi\right).
\end{align}
In the gap case, $u_*=0$, so that $\lambda_{\mathrm{c},2}$ is always positive, provided $\delta_2>0$. However, in the no-gap case, $u_*=1-\delta_1w_*$, so that
\begin{align*}
    \int_\mathbb{R} \delta_1u_*\bar{u}_*(\xi) &+\delta_2v_*(\xi) \bar{v}_*(\xi) \mathrm{d} \xi= \\ & \int_\mathbb{R}\left[ \frac{\delta_1(1-\delta_1w_*)}{c_*}\int^{\infty}_\xi e^{-\frac{1}{c_*}(-1+\delta_1w_*) (\xi-s)} v_{*,s}(s) \bar{v}_{*,s}(s)\mathrm{d}s +\delta_2v_*(\xi) v_{*,\xi}(\xi)e^{-\frac{c_*}{\kappa+\delta_1 w_*}\xi}\right] \mathrm{d} \xi.
\end{align*}
This term could be positive or negative depending on the relation between $\delta_1, \delta_2$, and the other system parameters $\rho, a, k, \delta_3$.

\section{Numerical simulations and discussion}
\label{S:NUM}
The formal geometric singular perturbation analysis of the preceding sections provides a framework by which we can understand the structure of $1$D bistable traveling tumor fronts in~\eqref{eq:GG}, and in particular uncovers geometric mechanisms which distinguish between qualitatively different cases: benign and malignant no-gap/gap tumors. Furthermore, when extending a $1$D profile as a straight planar interface in two spatial dimensions, the expression~\eqref{eq:lambda2cfinal} provides a stability criterion for long wavelength perturbations in the direction along the interface, assuming that the corresponding traveling wave is stable in one spatial dimension. 

\begin{figure}
\centering
\includegraphics[width=0.28\linewidth]{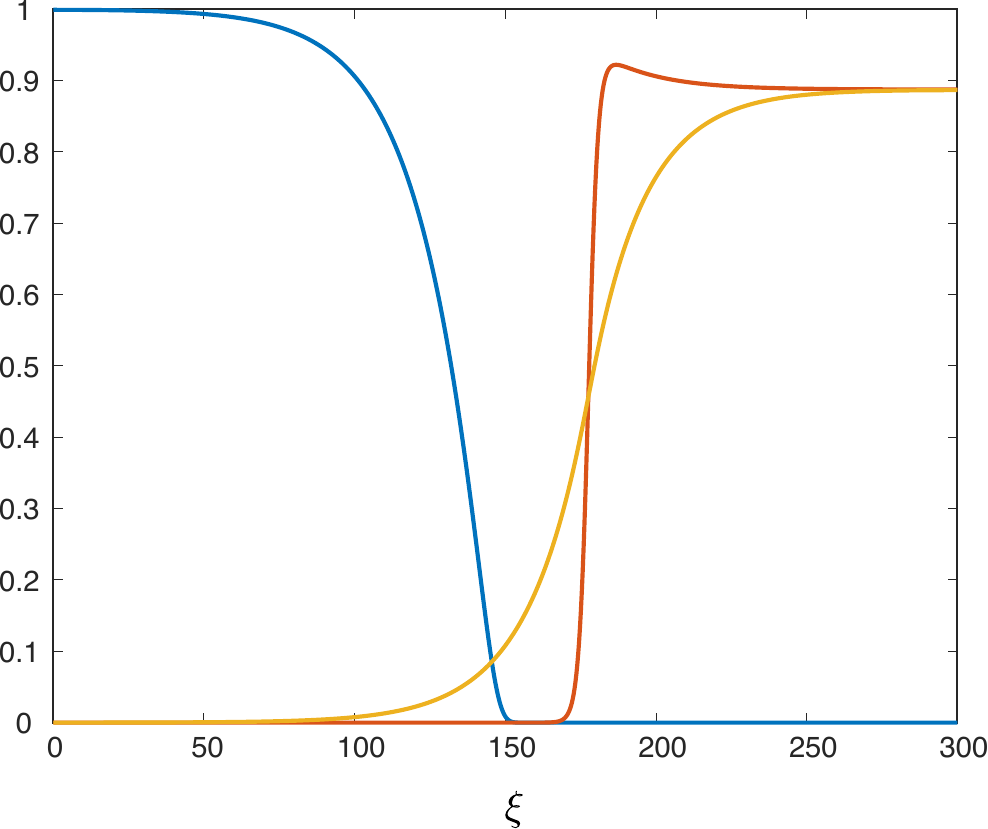}\hspace{0.025 \textwidth}
\includegraphics[width=0.3\linewidth]{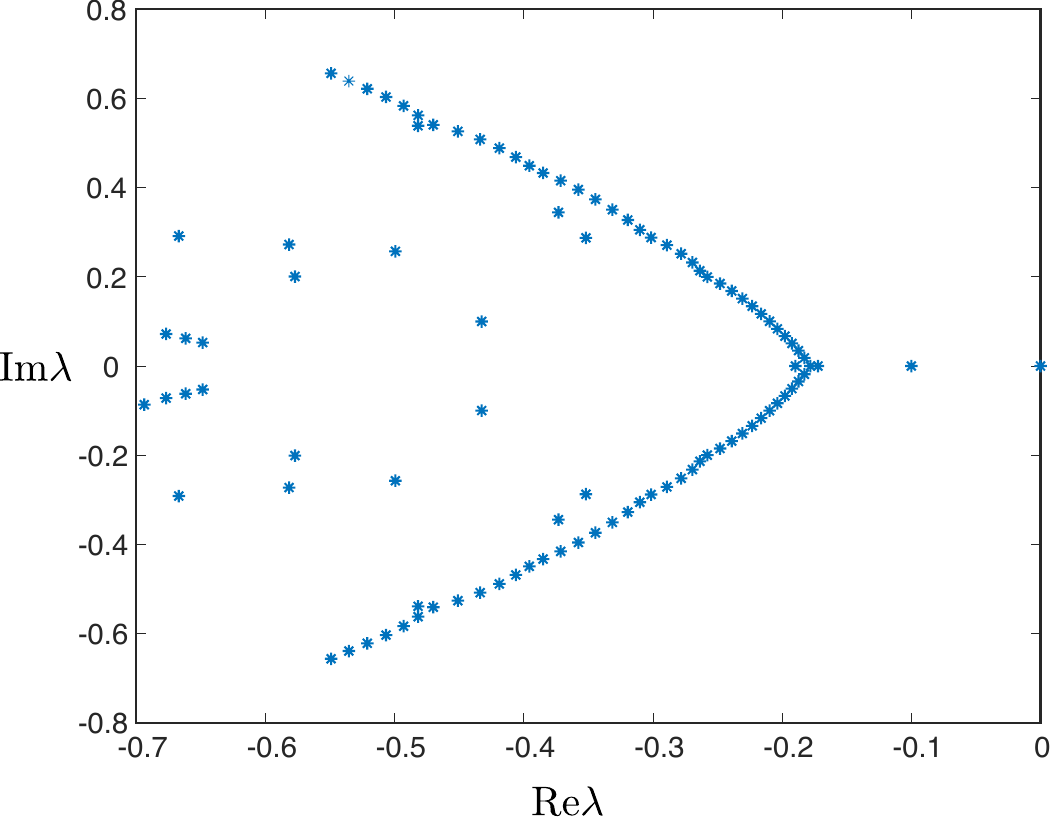}\hspace{0.025 \textwidth}
\includegraphics[width=0.3\linewidth]{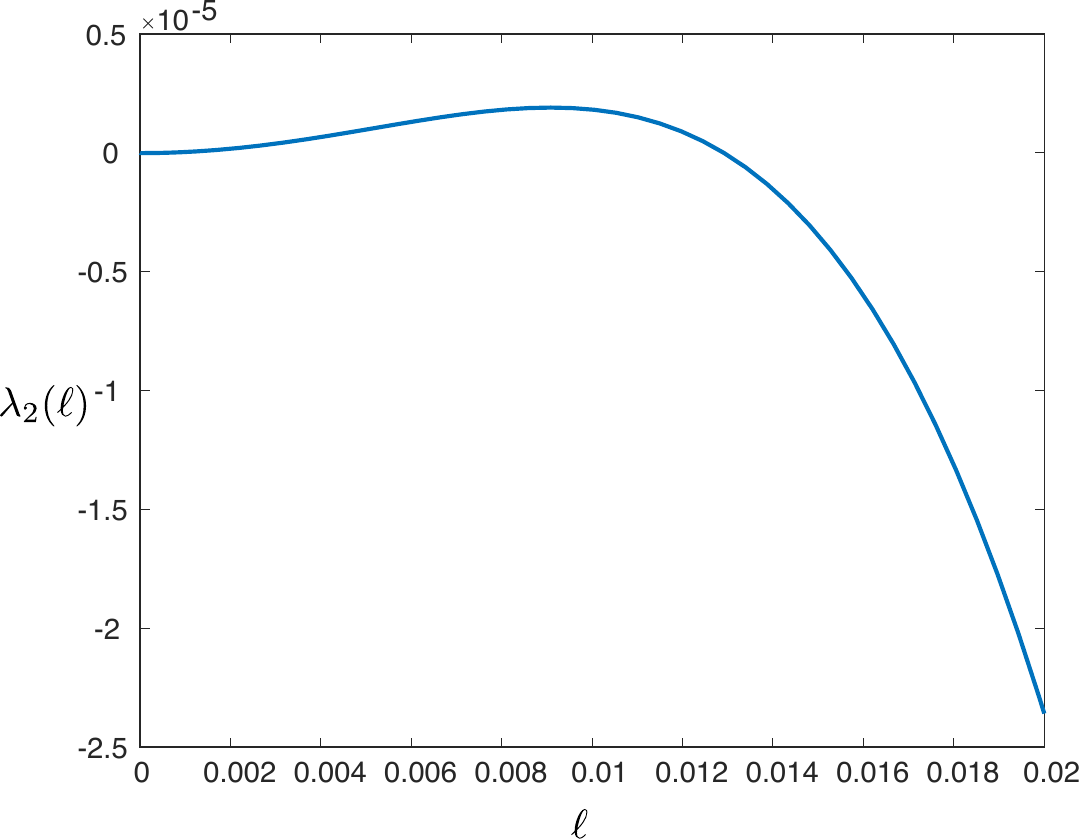}\\ 
\vspace{10pt}
\includegraphics[width=0.28\linewidth]{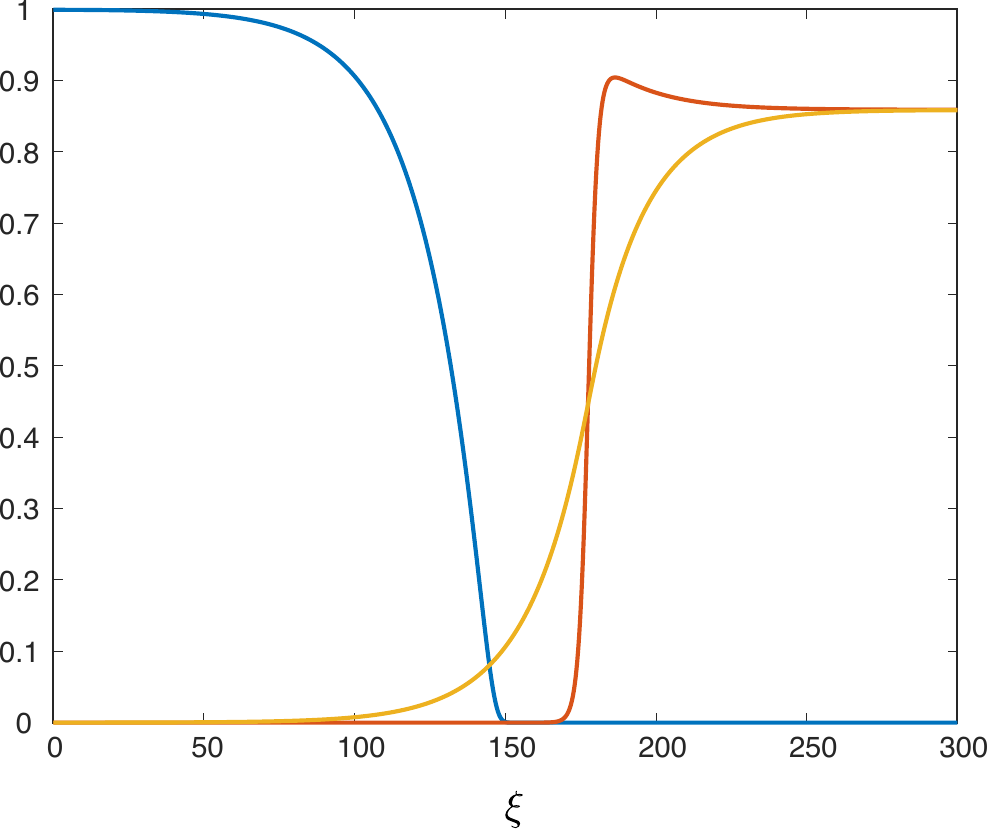}\hspace{0.025 \textwidth}
\includegraphics[width=0.3\linewidth]{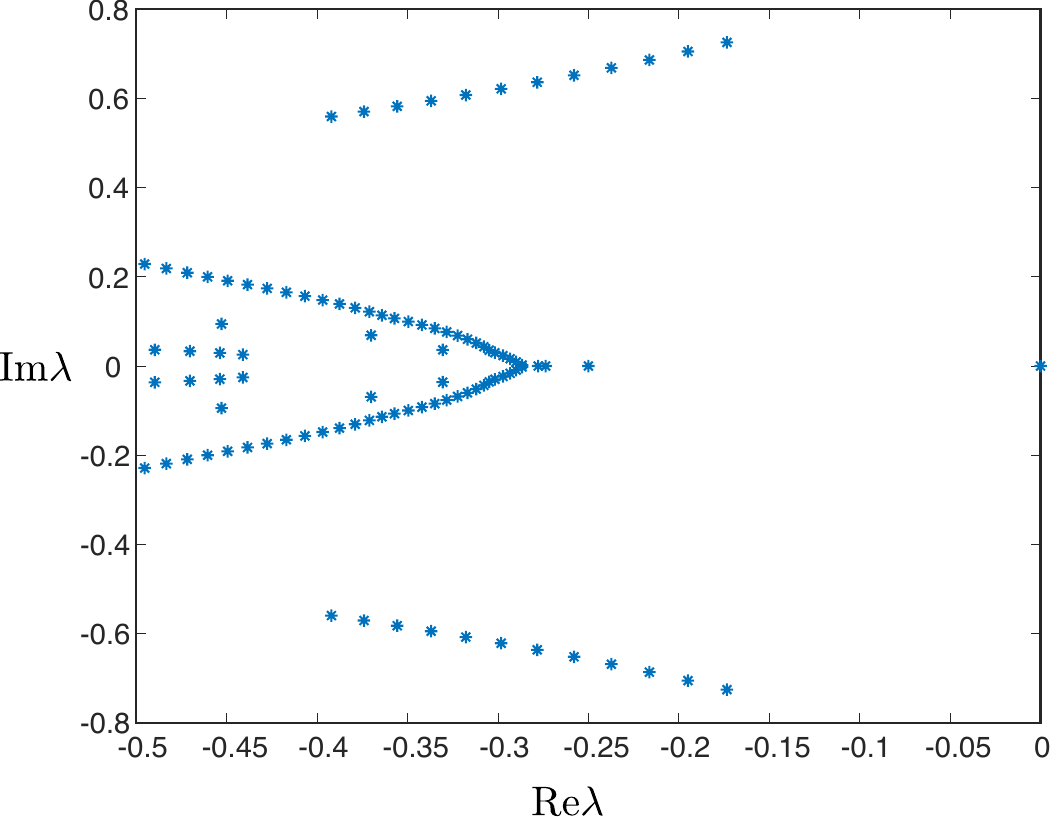}\hspace{0.025 \textwidth}
\includegraphics[width=0.3\linewidth]{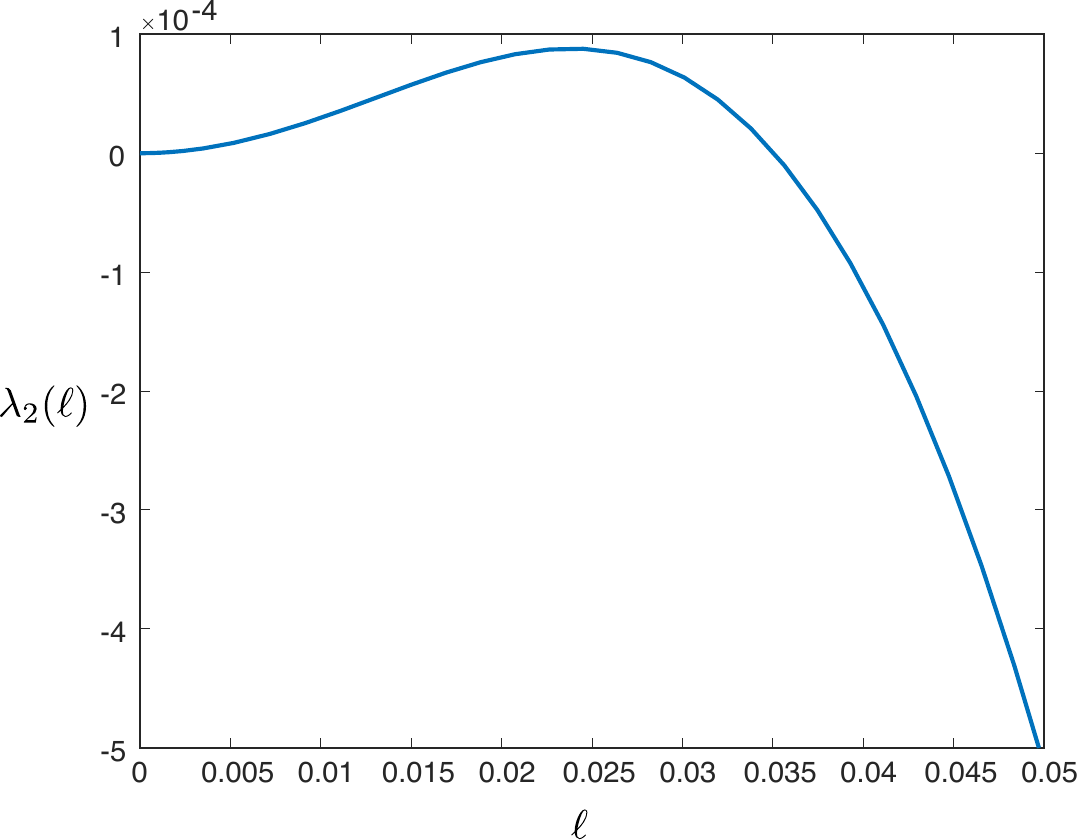}\\ \vspace{10pt}
\includegraphics[width=0.28\linewidth]{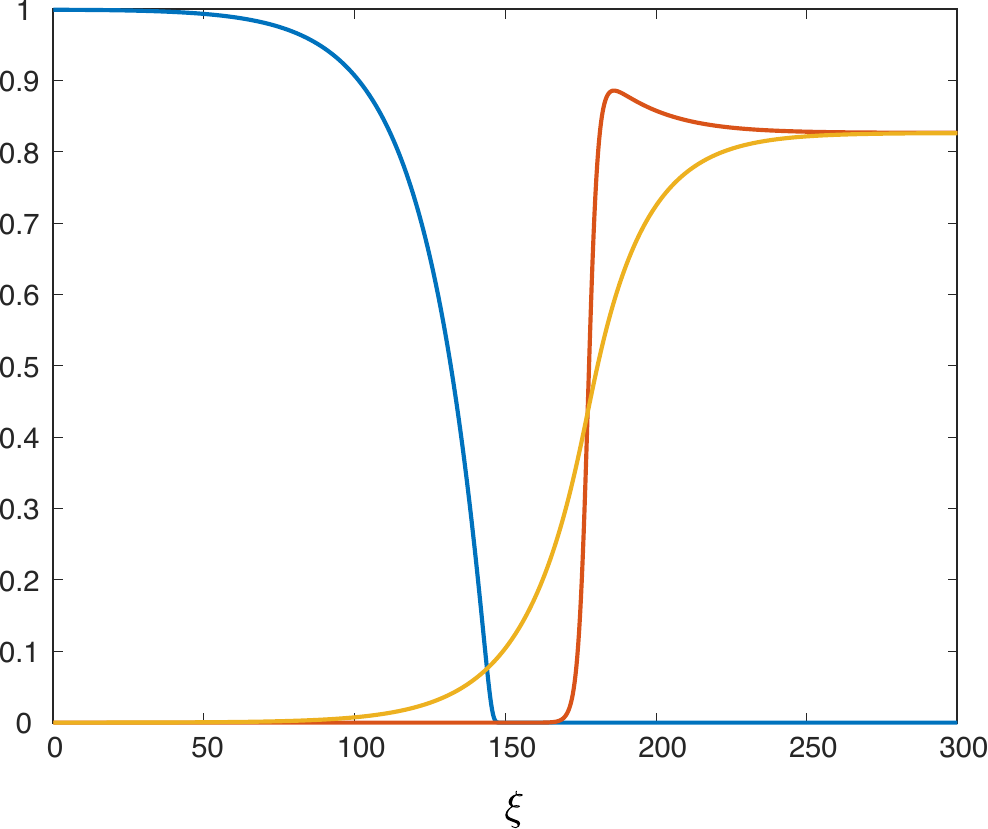}\hspace{0.025 \textwidth}
\includegraphics[width=0.3\linewidth]{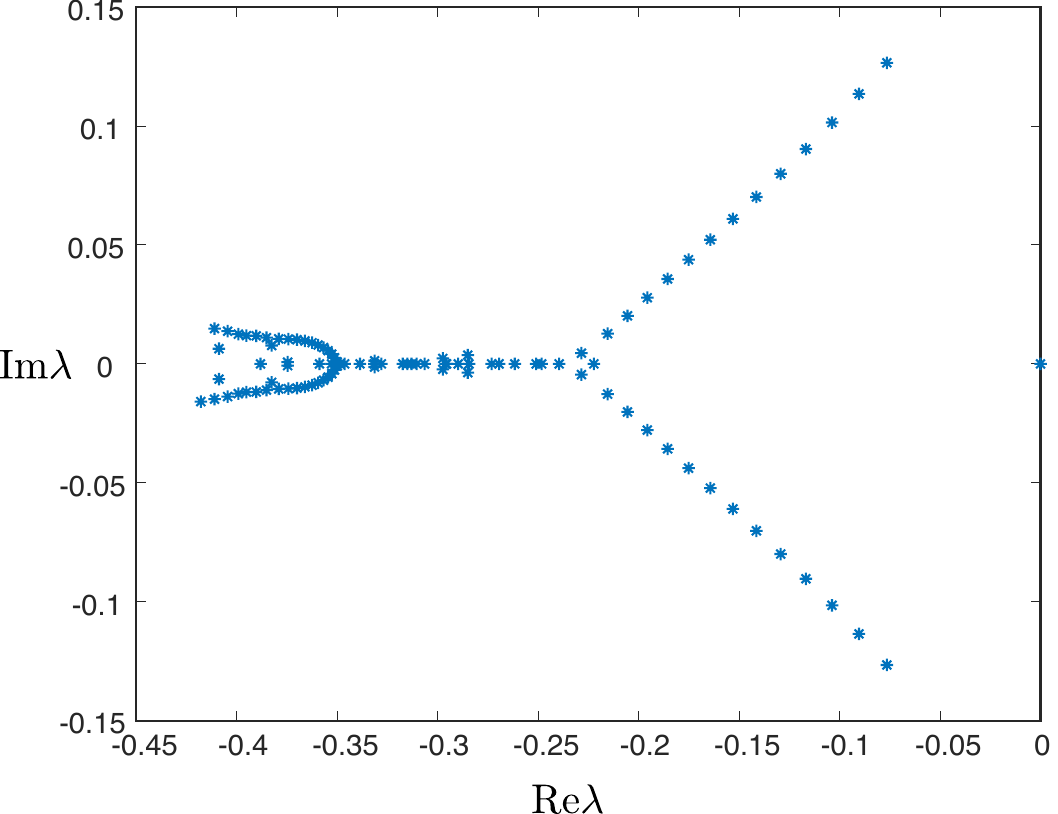}\hspace{0.025 \textwidth}
\includegraphics[width=0.3\linewidth]{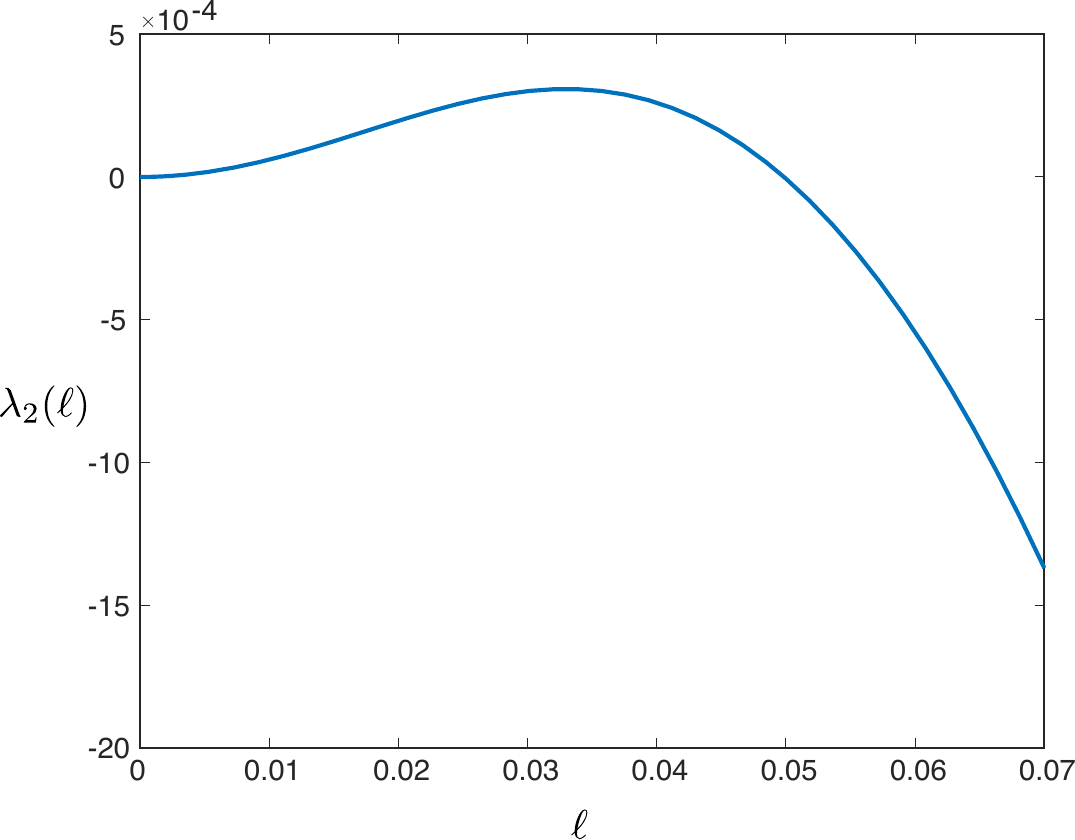}\\
\vspace{10pt} \includegraphics[width=0.28\linewidth]{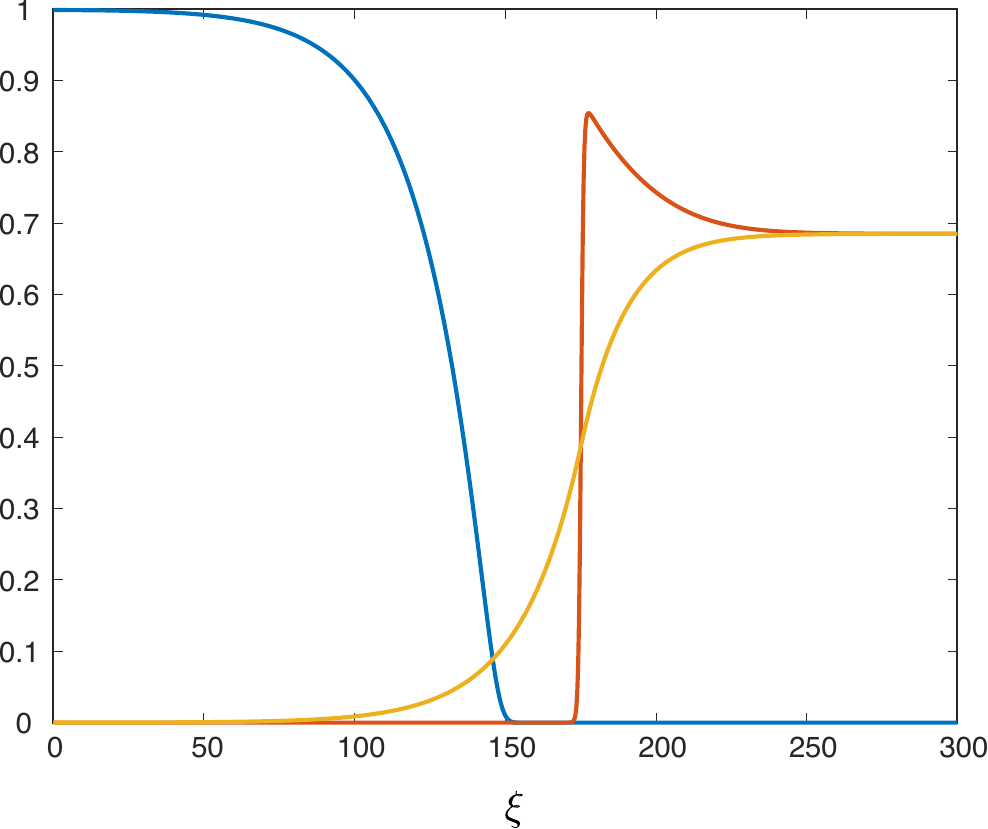}\hspace{0.025 \textwidth}
\includegraphics[width=0.3\linewidth]{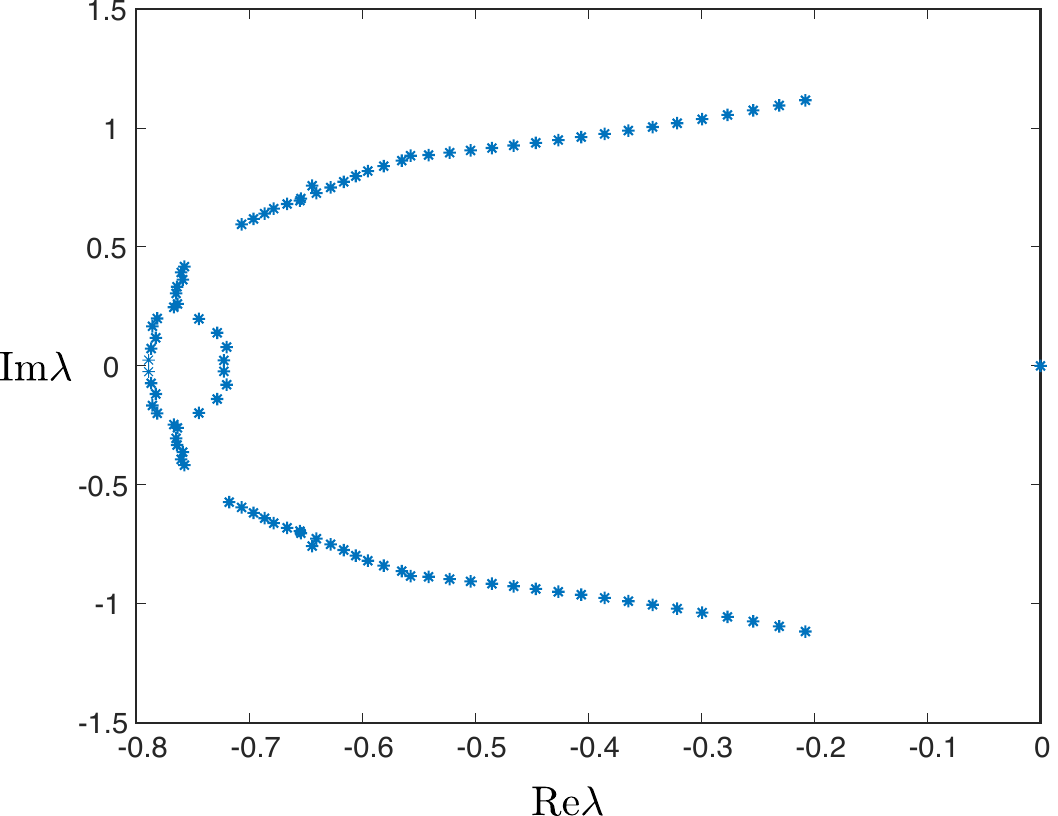}\hspace{0.025 \textwidth}
\includegraphics[width=0.3\linewidth]{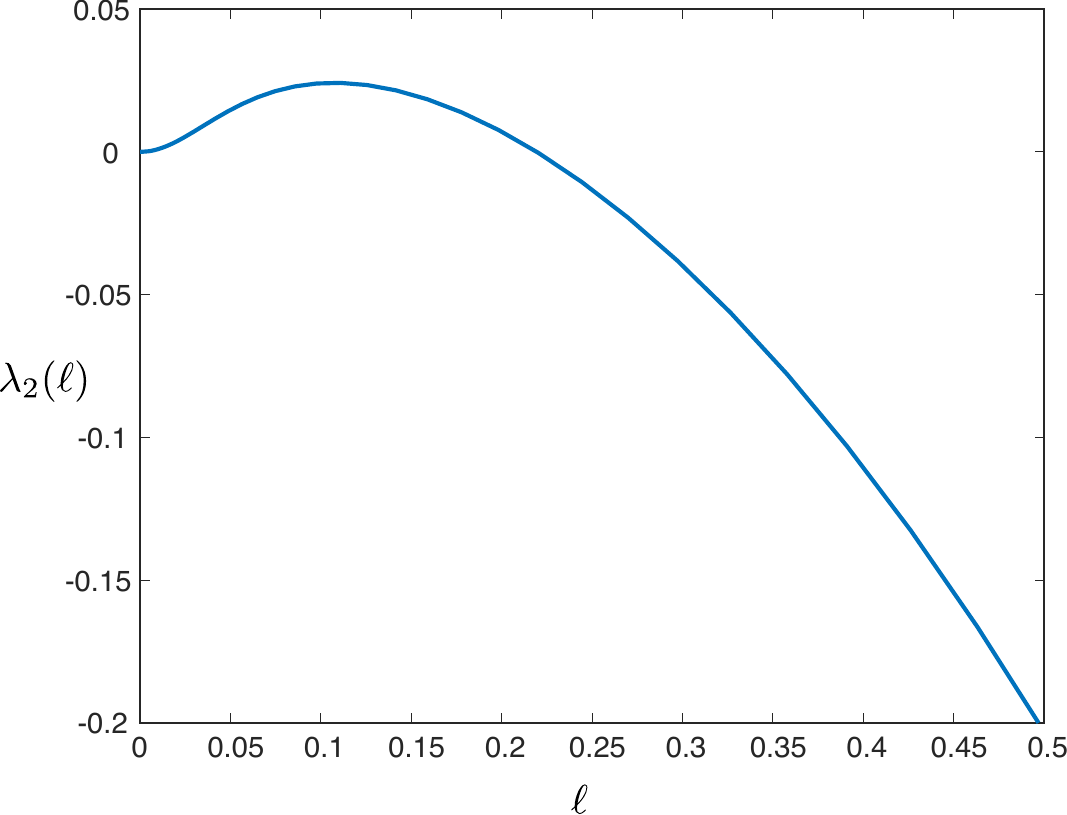}
\caption{$1$D traveling wave profiles obtained for the parameter values $(a,\kappa, \delta_1, \delta_2, \delta_3, \rho, \eps) = (0.1,0.1,12.5,0.1,70,1, 0.0063)$ (first row), 
$(a,\kappa, \delta_1, \delta_2, \delta_3, \rho, \eps) = (0.25,0.1,12.5,0.1,70,1, 0.0063)$ (second row), and $(a,\kappa, \delta_1, \delta_2, \delta_3, \rho, \eps) = (0.35,0.1,12.5,0.1,70,1, 0.0063)$ (third row), $(a,\kappa, \delta_1, \delta_2, \delta_3, \rho, \eps) = (0.25,0.05,11.5,3,1,15, 0.05)$ (fourth row). The $u,v,w$ profiles are plotted in blue, red, yellow, respectively. Profiles were obtained by solving the traveling wave equation~\eqref{eq:GG_twODE} in MATLAB. Also shown are the $1$D spectra, providing numerical evidence that all four solutions are $1$D-stable, as well as a continuation of the critical eigenvalue $\lambda_\mathrm{c}(\ell)$ for small, positive values of the wavenumber~$\ell$.
}
\label{f:profiles-spectra}
\end{figure}

While we can deduce from this expression and the arguments in~\S\ref{sec:lambda2c_asymp} that, for instance, tumor interfaces in the malignant gap case are always unstable in two spatial dimensions when $\eps$ is sufficiently small, in other regimes the sign of the integral expression~\eqref{eq:lambda2c_sgn} is less apparent due to its implicit dependence on various system parameters. However, we are able to explore different parameter regimes by numerically solving the traveling wave equation~\eqref{eq:GG_twODE}. Fig.~\ref{f:profiles-spectra} depicts numerically computed traveling waves profiles for four different sets of parameters. The $1$D spectra of the solutions are also shown, implying that these waves are stable in the longitudinal direction, so that when restricted to a one-dimensional domain, small perturbations of the wave decay in time. Also plotted is a continuation of the eigenvalue $\lambda_{\mathrm{c,2}}(\ell)$ for small values of the wave number $\ell$ in each case, indicating they are all unstable in two spatial dimensions; we note that each profile exhibits an acellular gap, and hence this matches our theoretical prediction that such solutions are unstable in two dimensions.

\begin{figure}
\centering
\includegraphics[width=0.39\linewidth]{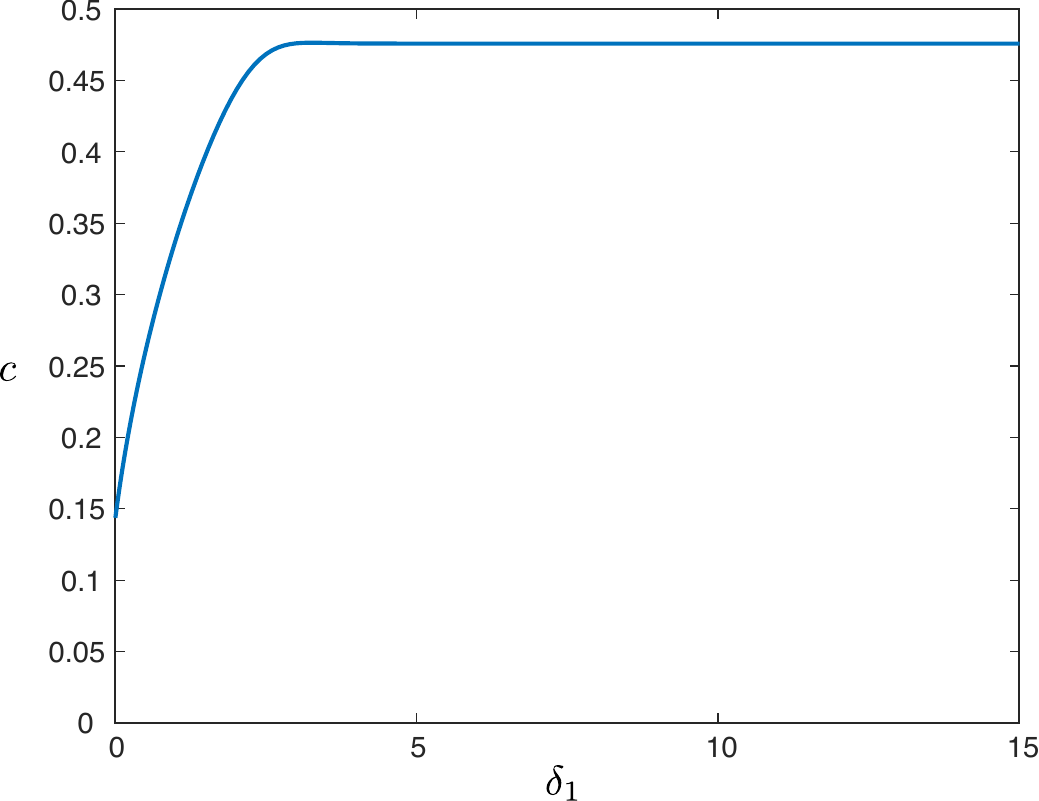}\hspace{0.1 \textwidth}
\includegraphics[width=0.4\linewidth]{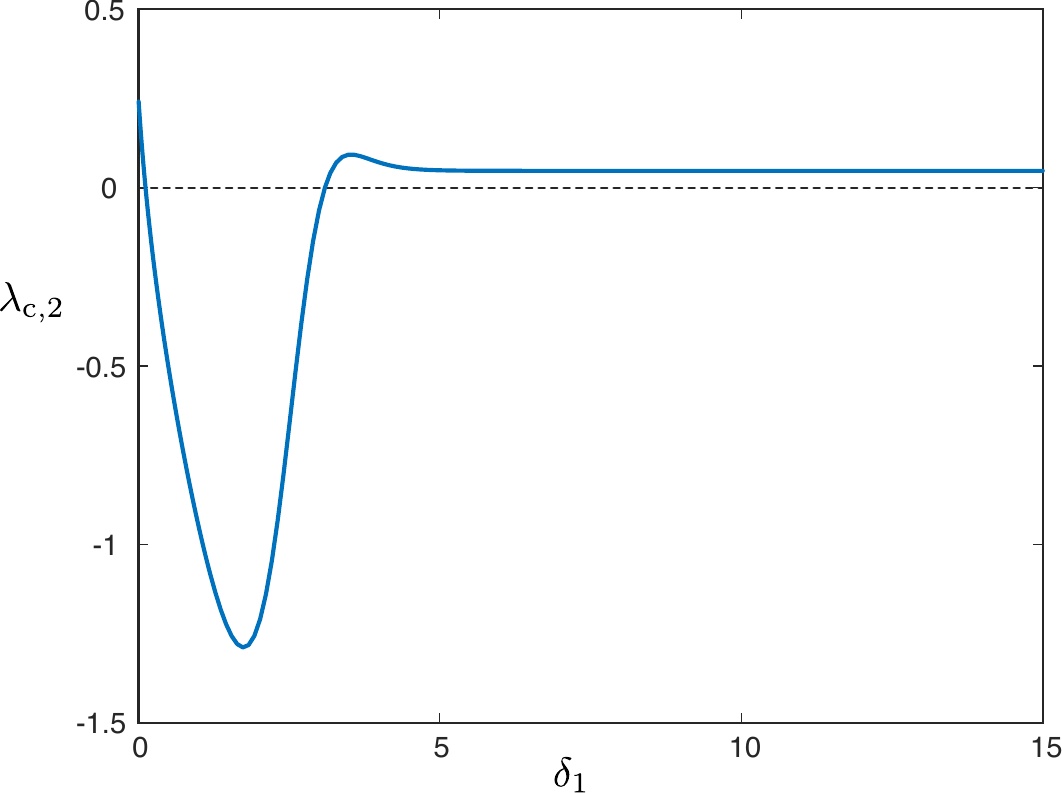}
\caption{Results of numerical continuation in AUTO07p\cite{doedel2007auto} for the parameter values $(a,\kappa, \delta_2,\delta_3, \rho, \eps) = (0.1,0.1,0.1,70,1.0,0.0063)$ for values of $\delta_1\in(0.05,15)$: wave speed $c$ versus $\delta_1$ (left), $\lambda_{\mathrm{c},2}$ versus $\delta_1$ (right) }
\label{f:delta1conthr}
\end{figure}

\begin{figure}
\centering
\includegraphics[width=0.39\linewidth]{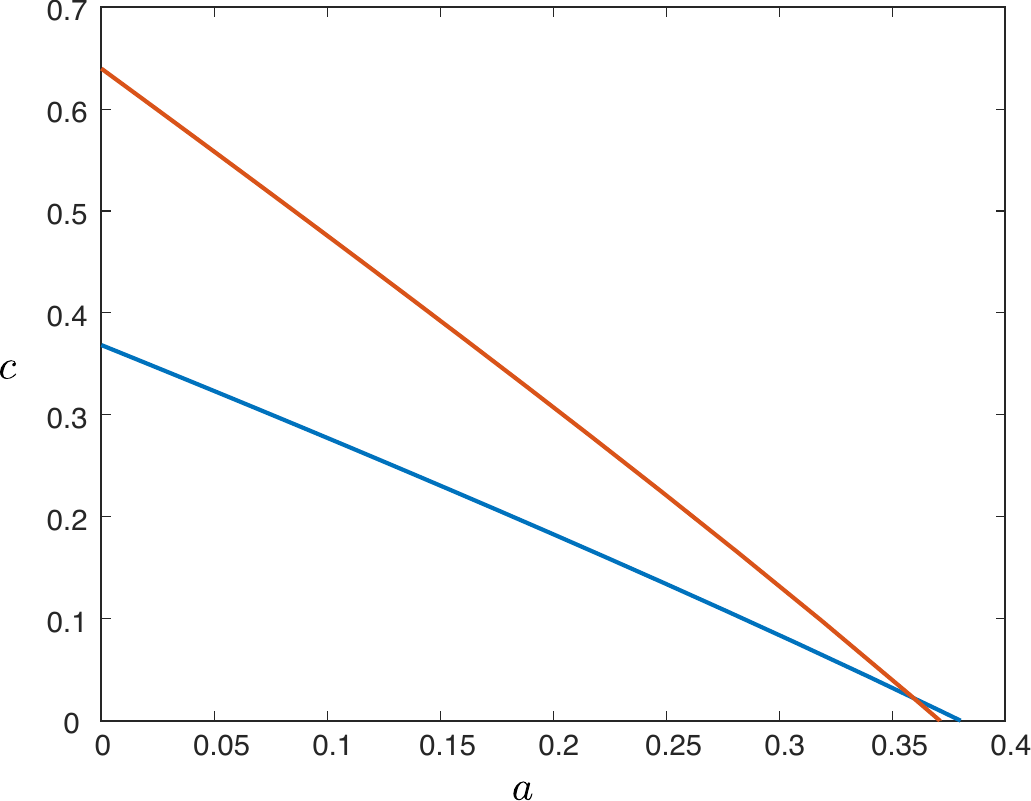}\hspace{0.1 \textwidth}
\includegraphics[width=0.4\linewidth]{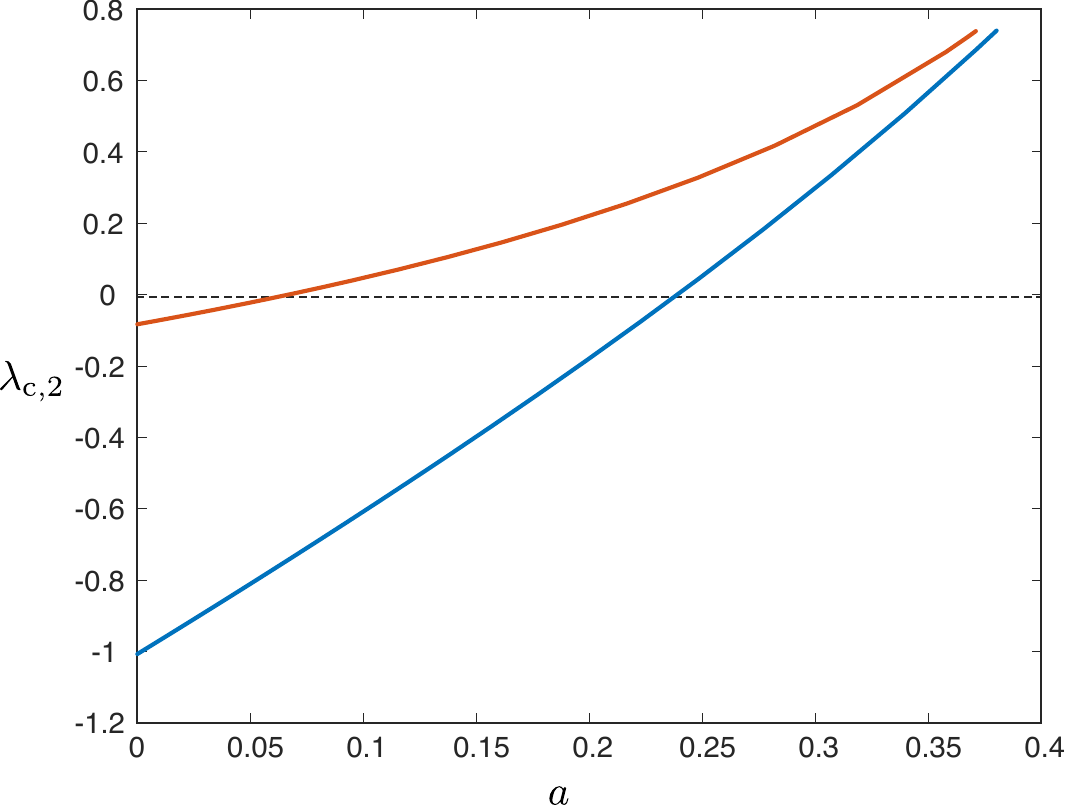}
\caption{Results of numerical continuation for the parameter values $(\kappa, \delta_2,\delta_3, \rho, \eps) = (0.1,0.1,70,1.0,0.0063)$ for a range of $a$-values for $\delta_1=0.6$ (blue) and $\delta_1=12.5$ (red): wave speed $c$ versus $a$ (left), $\lambda_{\mathrm{c},2}$ versus $a$ (right) }
\label{f:acont}
\end{figure}

Using the expression~\eqref{eq:lambda2c_expression}, we are also able to track this instability as a function of system parameters. Fig.~\ref{f:delta1conthr} depicts the results of numerical continuation in the parameter $\delta_1$ (the other parameters are the same as the wave in the first row of Fig.~\ref{f:profiles-spectra}). Note that while the waves can be $2$D stable for smaller values of $\delta_1$, for larger $\delta_1$ the instability appears, matching the prediction of the asymptotic expression~\eqref{eq:lambda2cfinal} that the interfaces are unstable in the malignant gap case -- the gap appears as $\delta_1$ increases; see Fig.~\ref{f:gapwidth}. We also note that the speed of the front also increases in $\delta_1$. In Fig.~\ref{f:acont}, we track the impact of the Allee effect on the coefficient $\lambda_{\mathrm{c},2}$. We see that in general, an increase in the Allee effect decreases the speed of the tumor, but also leads to the onset of the long wavelength instability. 

\begin{figure}
\centering
\includegraphics[width=0.4\linewidth]{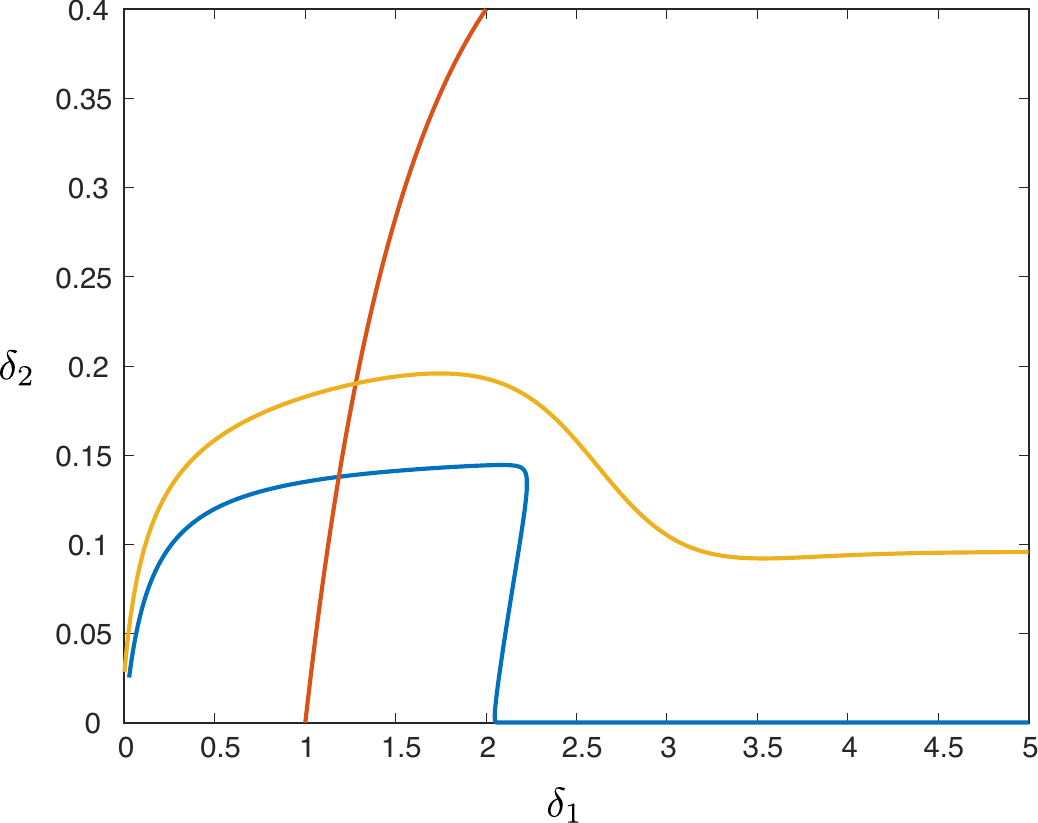}
\caption{Results of numerical continuation of the curve $\lambda_{\mathrm{c},2}=0$ for the parameter values $(a,\kappa, \delta_3, \rho) = (0.1,0.1,70,1.0)$ in the $(\delta_1,\delta_2)$-plane for $\eps=0.0063$ (yellow) and $\eps=10^{-5}$ (blue). Below each curve we have $\lambda_{\mathrm{c},2}<0$, while $\lambda_{\mathrm{c},2}>0$ above. Plotted in red is the malignant/benign boundary $\delta_1V_+=1$.  }
\label{f:boundary}
\end{figure}

By solving numerically for zeros of the expression~\eqref{eq:lambda2c_expression}, we can also track the stability boundary in parameter space. Fig.~\ref{f:boundary} depicts the stability boundary in $(\delta_1,\delta_2)$-space for two different values of $\eps$. Also shown is the curve which forms the boundary in parameter space between the benign and malignant tumors, which shows that both benign/malignant tumors can be stable/unstable and that in general, the interfaces are more likely to be unstable for larger values of $\delta_1, \delta_2$. 

\begin{figure}
\centering
\includegraphics[width=0.22\linewidth]{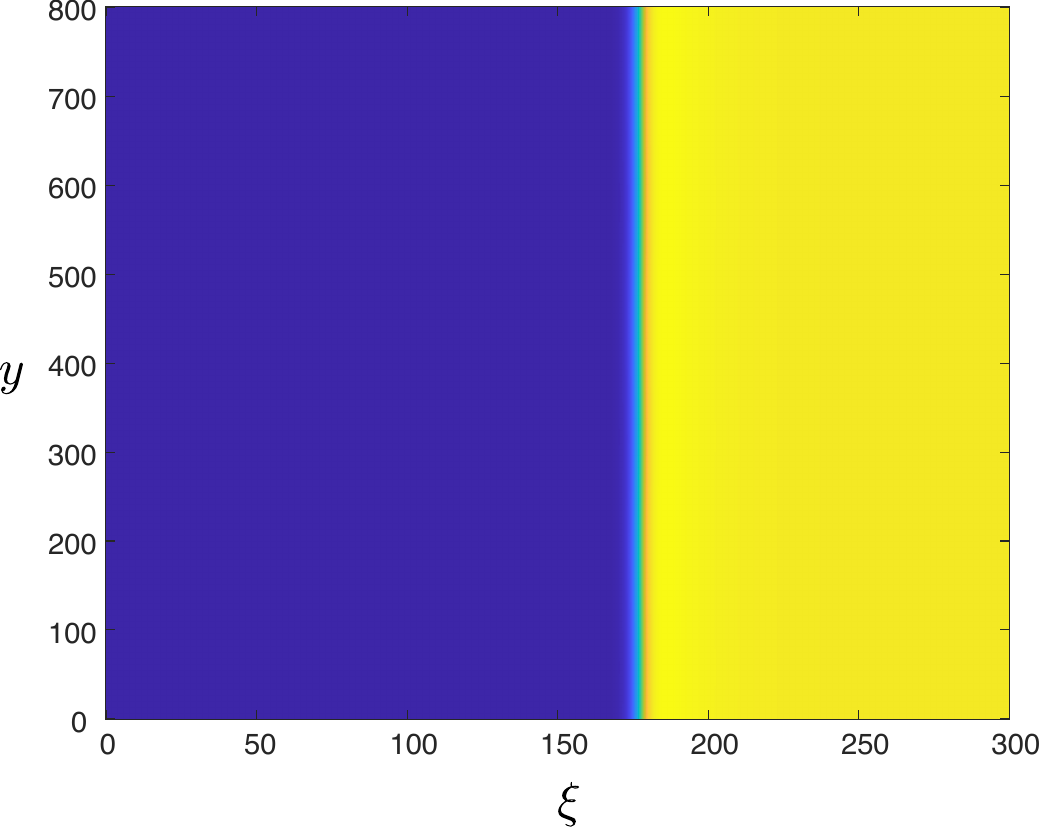}\hspace{0.025 \textwidth}
\includegraphics[width=0.22\linewidth]{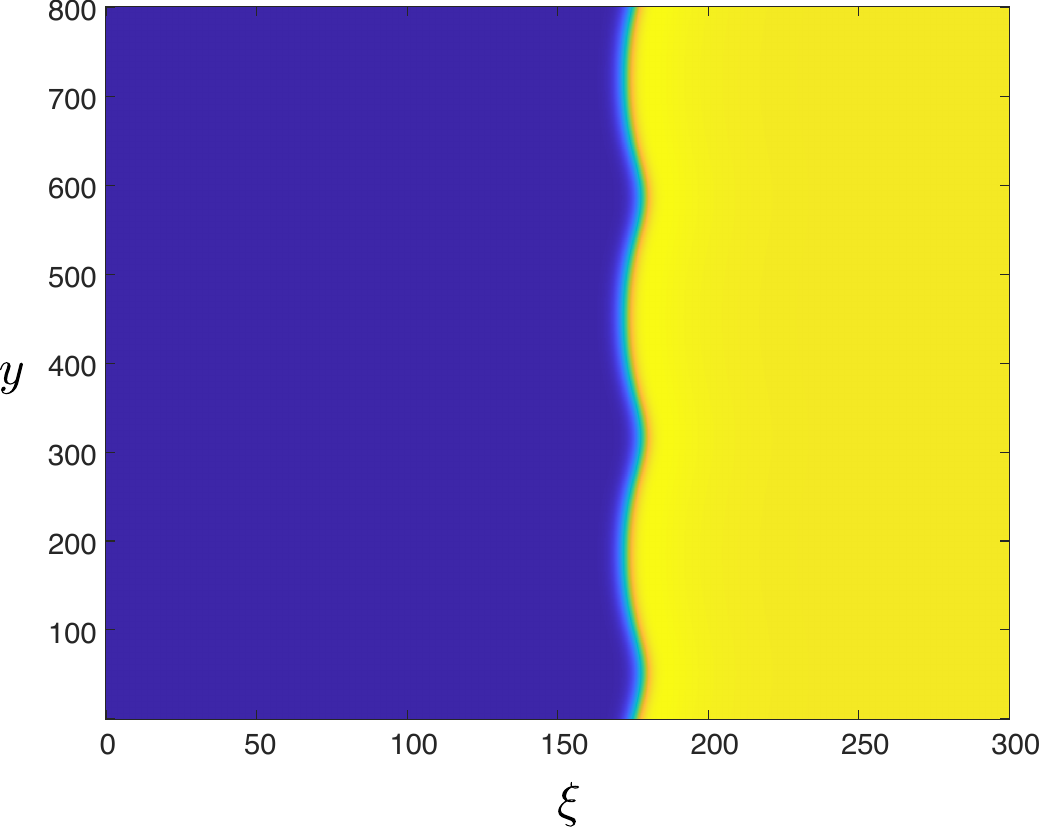}\hspace{0.025 \textwidth}
\includegraphics[width=0.22\linewidth]{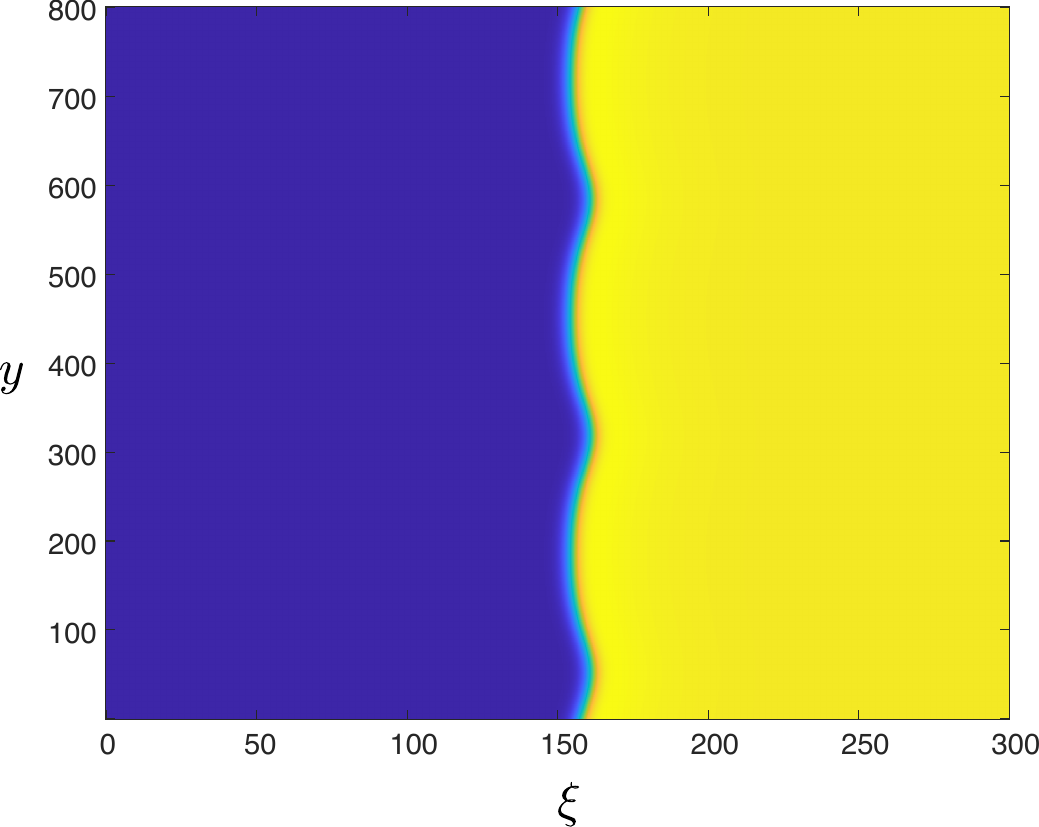}\hspace{0.025 \textwidth}
\includegraphics[width=0.22\linewidth]{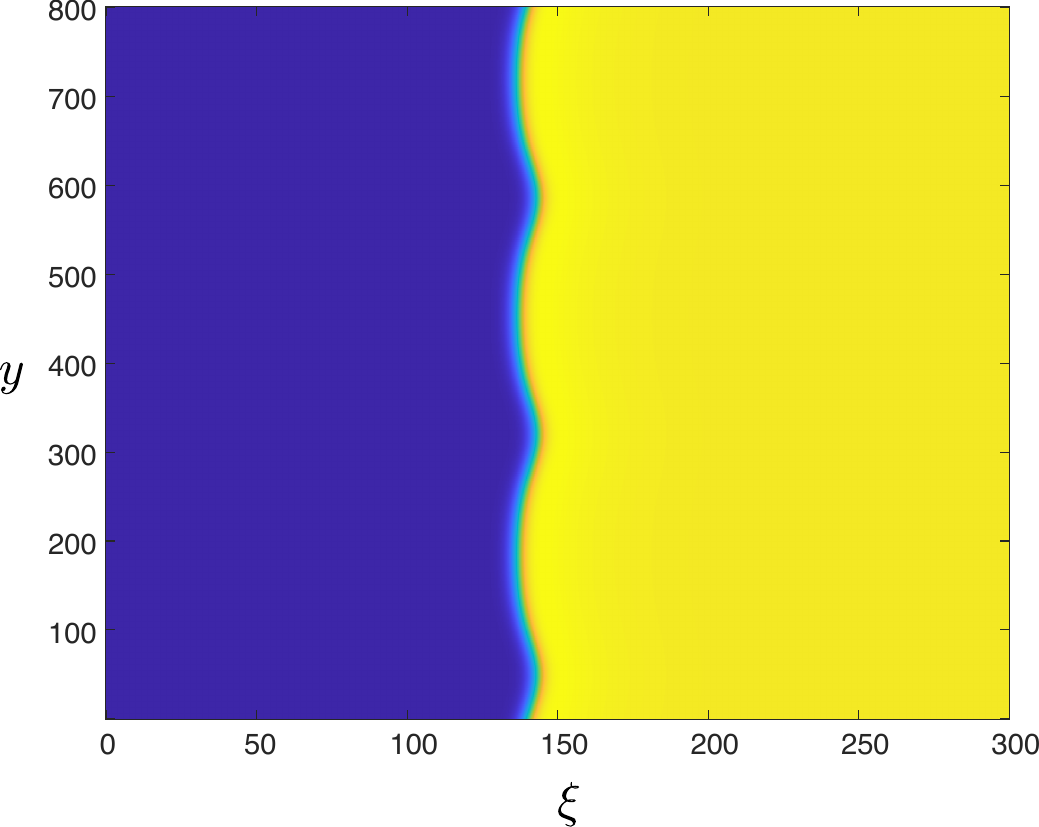}\\
\includegraphics[width=0.22\linewidth]{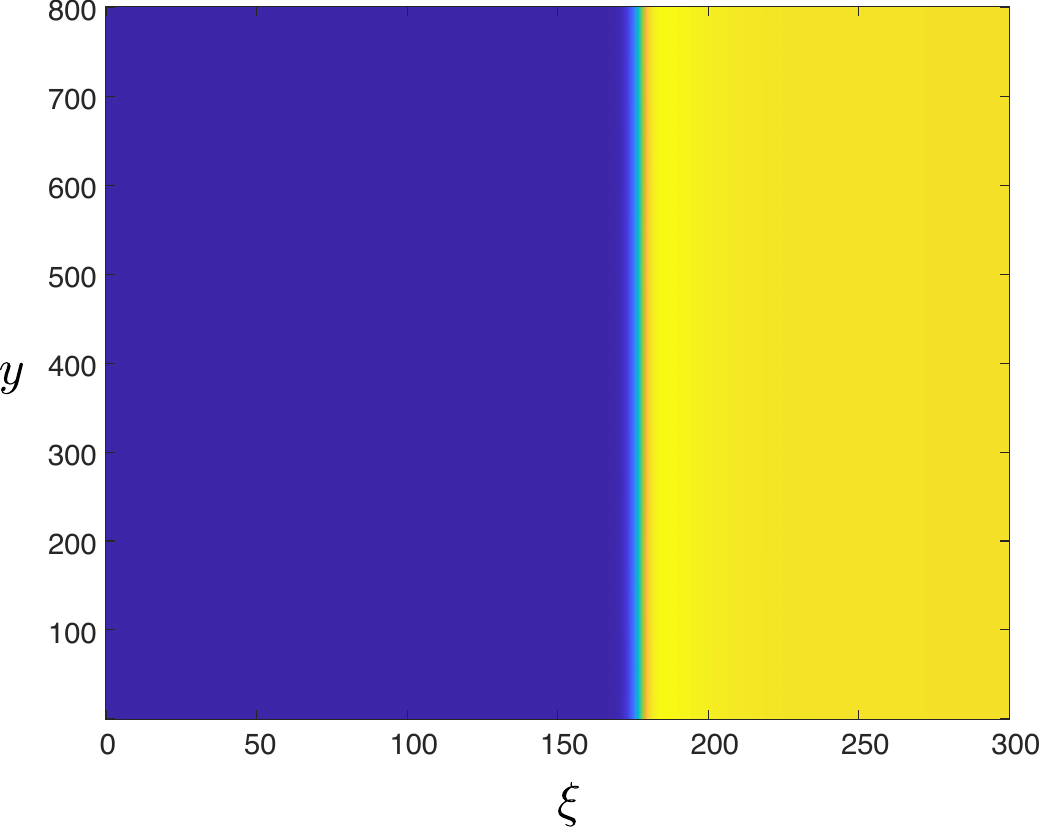}\hspace{0.025 \textwidth}
\includegraphics[width=0.22\linewidth]{Figures/a35_t_20000.pdf}\hspace{0.025 \textwidth}
\includegraphics[width=0.22\linewidth]{Figures/a35_t_30000.pdf}\hspace{0.025 \textwidth}
\includegraphics[width=0.22\linewidth]{Figures/a35_t_40000.pdf}\\
\includegraphics[width=0.22\linewidth]{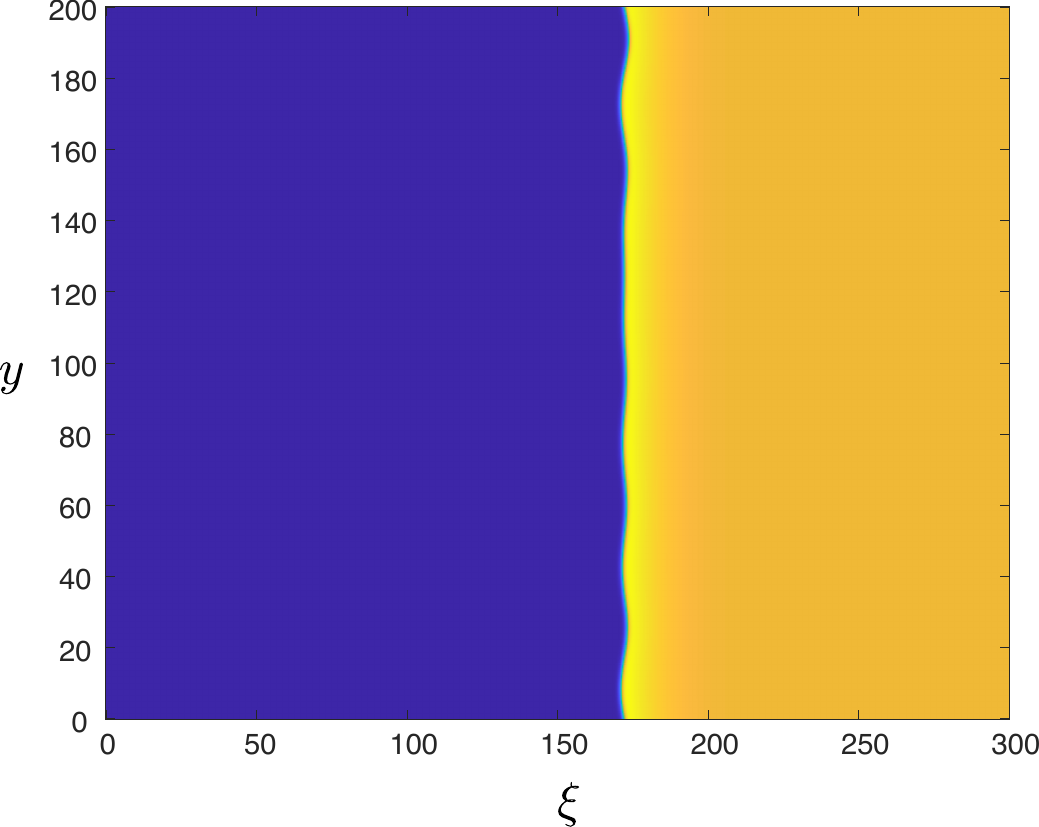}\hspace{0.025 \textwidth}
\includegraphics[width=0.22\linewidth]{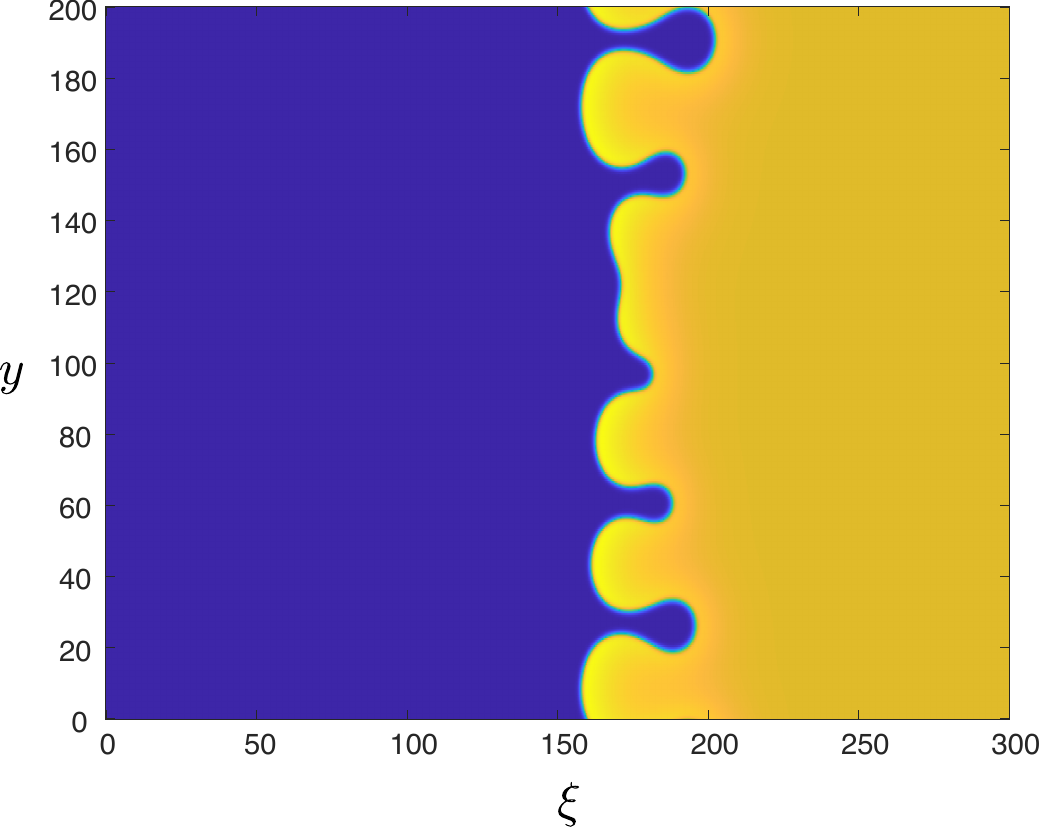}\hspace{0.025 \textwidth}
\includegraphics[width=0.22\linewidth]{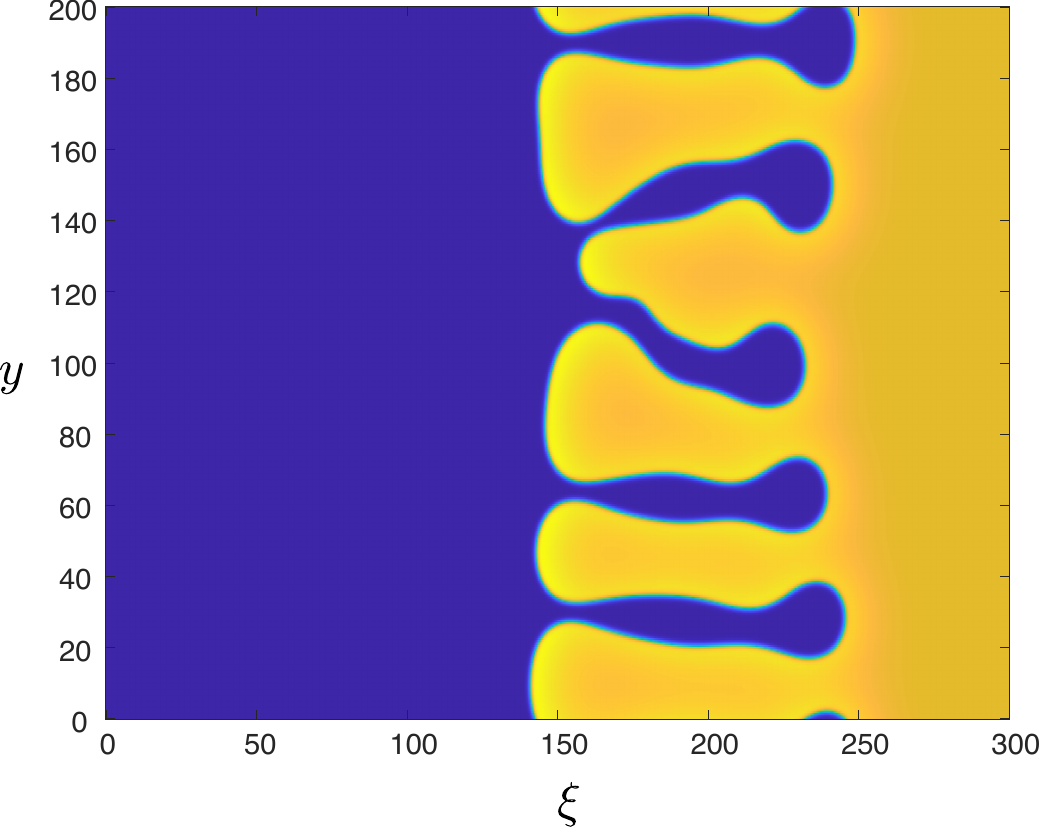}\hspace{0.025 \textwidth}
\includegraphics[width=0.22\linewidth]{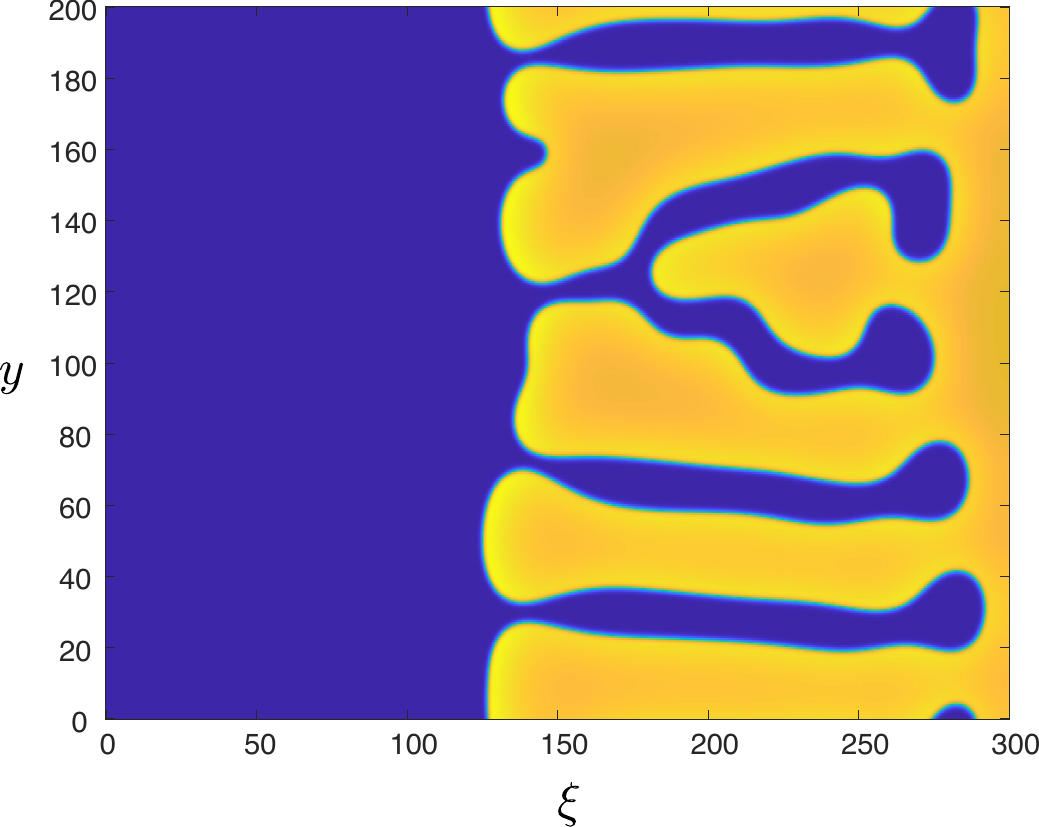}
\caption{Results of direct numerical simulations exhibiting different manifestations of the long wavelength instability. The simulations were performed in a co-moving frame corresponding to the wave speed of the initial front profile (so that in the absence of any instability, the fronts would appear stationary), with Neumann boundary conditions in $\xi$, and periodic boundary conditions in $y$. Finite differences were used for spatial discretization, and  MATLAB's ode15s routine was used for time stepping. The corresponding $v$-profile is depicted, where yellow indicates high density of tumor cells, and blue depicts low density of tumor cells. 
(First row) Simulation for the parameter values $(a,\kappa, \delta_1, \delta_2, \delta_3, \rho, \eps) = (0.25,0.1,12.5,0.1,70.0,1.0, 0.0063)$ at the times $t=5000,10000,15000,20000$ from left to right, with wave speed $c=0.2211$. (Second row) Simulation for the parameter values $(a,\kappa, \delta_1, \delta_2, \delta_3, \rho, \eps) = (0.35,0.1,12.5,0.1,70.0,1.0, 0.0063)$ at the times $t=10000,20000,30000,40000$ from left to right, with wave speed $c=0.0401$. (Third row) Simulation for the parameter values $(a,\kappa, \delta_1, \delta_2, \delta_3, \rho, \eps) = (0.25,0.05,11.5,3,1,15, 0.05)$ at the times $t=250,350,450,550$ from left to right, with wave speed $c=0.3296$.  In each of the three simulations, initial conditions were constructed by trivially extending the $1$D-stable profiles from Fig.~\ref{f:profiles-spectra} (rows $2$ through $4$) in the $y$-direction and adding a small amount of positive noise. }
\label{f:interface}
\end{figure}

We emphasize that our analysis concerns spectral stability of the traveling fronts, and while these spectral computations may indicate instability, it is not easy to predict the resulting nonlinear behavior of the interface. Nevertheless, direct numerical simulations suggest that planar interfaces are in fact nonlinearly unstable, and the dynamics may lead to complex patterns at the tumor interface. Fig.~\ref{f:interface} depicts the results of direct numerical simulations for three different choices of the parameters $(a,\eps)$, initialized with $1$D-stable traveling wave profiles taken from rows two through four of Fig.~\ref{f:profiles-spectra}. Different manifestations of the long wavelength instability appear in each case: in the first simulation the interface develops cusps which remain bounded as the propagation speed of the advancing interface appears to increase, while in the second, the interface breaks up into growing finger patterns; the parameter values for these two simulations were taken from~\cite{GG} (with the exception of $(a,\delta_2,\kappa)$, which were not present in~\cite{GG}). The third and final simulation exhibits more elaborate growing finger patterns, which develop on a faster timescale. Note that the time taken for the instability to develop varies based on the magnitude of the coefficient $\lambda_{\mathrm{c},2}$; see Fig.~\ref{f:profiles-spectra}.

From the mathematical point of view, there are a number of further research directions that are directly in line with the present work. More general models also include tumor cell density $V$ in the nonlinear diffusion term of the evolution equation for $V$ \eqref{eq:GG} and may also contain (nonlinear) diffusion on the $U$-equation. The methods developed here are believed to be sufficiently flexible to investigate more general equations. Within the current framework, several open questions remain. In particular, the basic assumption underlying the two-dimensional stability analysis here (that is numerically validated in several cases, see Fig. \ref{f:profiles-spectra}) is the stability of the front with respect to one-dimensional (longitudinal) perturbations -- a rigorous verification of this assumption is the subject of future work. In addition, our (in)stability results are purely at the spectral level, and a challenging direction for future work concerns the study of the nonlinear manifestation of the long wavelength instability studied here.

Our mathematical analysis and numerical simulations naturally lead to a number of biological insights, some of which we now highlight. We showed that the speed of invasion, and the probability of an acellular gap being formed, both increase as $\delta_1$ increases. As this parameter is a measure of the negative impact of lactic acid on normal cells, this result makes intuitive sense. Furthermore, the stronger the Allee effect on the tumor cells, the slower the invasion speed of these cells, a result that also aligns with our intuition. Our study of the front instabilities showed that in such cases (strong Allee effect), the front undergoes a bifurcation that initiates the formation of growing `fingers', resulting in an irregular morphology, while in the case of a weak Allee effect, we predict that the waves will move faster and the initial bifurcation has a milder effect and drives the development of (moving) `cusps' in the front morpholgy, similar to that observed in the experimental results of Gatenby and Gawslinski \cite{GG}. 

Note that in our model, the instability of the invasion tumor cell front is an emergent property of the system which does not require phenotypic heterogeneity within the tumor cell population or spatial heterogeneity within the system (these are often assumed to be the cases for the break-up of an invading front). However, in reality, tumors typically are very heterogeneous and the extracellular matrix (ECM) through which they move is not spatially homogeneous. Possible future avenues of research would be to extend our present work to include these complications. Again, there is a strong similarity between cancer research and ecology (cf. \cite{korolev2014turning}), in which understanding the impact of spatial heterogeneities also is a central issue (cf. \cite{bastiaansen2019dynamics} and the references therein). In fact, through the link with ecology, another promising (and challenging) line of research emerges: like in ecology, interfaces between different states -- bare soil/vegetated or normal/tumorous -- are typically curved, not flat (as assumed here). Thus, it is necessary to extend the current approach to curved interfaces -- see \cite{byrnes2023large} for some first steps in that direction in the ecological setting. Note that there also is an important distinction between ecosystem and tumor interfaces: the former typically are within two-dimensional domains (the surface of a terrain), while the latter are intrinsically three-dimensional: in combination with the local curvatures, this additional freedom may for instance have a significant impact on the nature of the protrusions initiated by the fingering instability.    

From the modeling point of view, there are also various promising future research directions. For example, \cite{Strobl20} presents a model of cancer cell invasion in which there are two different tumor cell phenotypes, one producing lactic acid and the other producing proteins that degrade the ECM, considered as a barrier to invasion. They show how different phenotypic spatial structures arise in the invading front depending on the inter-cellular competition dynamics between the two cancer cell phenotypes. More recently, Crossley et al. \cite{crossley2024phenotypic} present a model based on volume-filling, in which one cancer cell phenotype proliferates, while the other degrades the matrix (an example of the well-known “go-or-grow” hypothesis). Analysis of this model also gives rise to different phenotypically structured invading fronts, depending on parameter values. Both of these studies were carried out in one spatial dimension and it would be interesting to see what structures, in both physical and phenotypic space, they exhibit when considered in two (or three) spatial dimensions. A further modelling complication to address is that, while these models consider the ECM as a barrier to cancer cell invasion, the ECM also enables cell invasion through providing a “scaffold” to which cells can attach and move. Understanding in detail how this dual property of the ECM affects tumor invasion is an open question.

\paragraph{Acknowledgments:} PC, DL, EO, and PY were partially supported by the National Science Foundation through grant DMS-2105816 and acknowledge the hospitality of the Mathematical Institute at Leiden University where much of the research for this project was carried out. The research of AD is supported by the ERC-Synergy project RESILIENCE (101071417). PvH acknowledges support by the Australian Research Council (DP190102545). PKM would like to thank the NSF-Simons Center for
Multiscale Cell Fate at the University of California, Irvine, for partial support.  The authors thank Robert Vink (Laboratory of Pathology East Netherlands, Hengelo, The Netherlands) for his insightful and helpful comments about the nature of observed tumors. 

\paragraph{Data availability:} Simulation data and codes generated during the current study are available from the corresponding author on reasonable request.

\appendix

\section{Stability of steady states}\label{app:steadystates}
We consider
\begin{align*}
\begin{split}
U_\tau &=F(U,W) \,,\\
V_\tau  &= G(U,V,W)+\nabla\cdot((1+\kappa-U)\nabla V)\,,\\
W_\tau &= H(V,W)+\frac{1}{\eps^2}\Delta W  \,,
\end{split}
\end{align*}
where $F,G,H$ are as in~\eqref{eq:FGHdef}. Linearizing about a steady state $(U,V,W) = (U_0,V_0,W_0)+(\bar{u}, \bar{v}, \bar{w})e^{i(\ell_1 x+\ell_2 y)+\lambda t}$, we obtain the linearized system
\begin{align}
\lambda \begin{pmatrix} \bar{u} \\ \bar{v} \\ \bar{w}\end{pmatrix}= \begin{pmatrix}F_{u} & 0 & F_w \\ 0 & -(1+\kappa-U_0)(\ell_1^2+\ell_2^2)+ G_{v}  & G_w \\ 0 & H_{v} & -\frac{1}{\eps^2}(\ell_1^2+\ell_2^2) +H_w   \end{pmatrix}\begin{pmatrix} \bar{u} \\ \bar{v} \\ \bar{w}\end{pmatrix}\,,
\end{align}
where $F_u: = \tfrac{\partial F}{\partial u}(U_0,V_0, W_0)$, etc. The steady state is stable if all eigenvalues satisfy $\mathrm{Re}(\lambda)<0$ for all $\ell_1, \ell_2\in \mathbb{R}$. From this we see that a necessary condition is $F_u<0$, from which we immediately obtain that the steady state $P_1$ is always unstable, while $P_3^+$ is unstable if $\delta_1 V^+>1$ and $P_4^+$ is unstable if $\delta_1 V^+<1$. 

The stability of the remaining steady states is determined by the lower right $2 \times 2$ block eigenvalue problem
\begin{align}
\lambda \begin{pmatrix} \bar{v} \\ \bar{w}\end{pmatrix}= \begin{pmatrix}-(1+\kappa-U_0)(\ell_1^2+\ell_2^2)+ G_{v}  & G_w \\ H_{v} & -\frac{1}{\eps^2}(\ell_1^2+\ell_2^2) +H_w   \end{pmatrix}\begin{pmatrix}  \bar{v} \\ \bar{w}\end{pmatrix}\,,
\end{align}
which, following~\cite[\S2.1]{CDLOR}, are stable for $0<\eps\ll1$, provided the conditions 
\begin{align}
  G_v<0, \qquad  G_v+H_w<0, \qquad G_vH_w-G_wH_v>0\,,
\end{align}
are satisfied. Employing these conditions, a short computation shows that $P_2$ is stable, $P_3^+$ is stable if $\delta_1 V^+<1$ and $P_4^+$ is stable if $\delta_1 V^+>1$, and the remaining steady states are unstable.

\bibliographystyle{abbrv}
\bibliography{my_bib}

\end{document}